\input amstex
\input epsf
\documentstyle{amsppt}
\magnification=1200
\NoBlackBoxes
\loadmsbm
\scriptscriptfont0=\sevenrm
\scriptscriptfont1=\seveni
\scriptscriptfont2=\sevensy
%\pagewidth{143mm}
%\pageheight{205mm}

\define\ap {{\text {\rm '}}}
\hcorrection{.25 in}
\vcorrection{ .5 in}
\advance\vsize-.75in.
\ifx\undefined\rom
\define\rom#1{{\rm #1}}
\fi
%  Macro for current address.
\ifx\undefined\curraddr
\def\curraddr#1\endcurraddr{\address {\it Current address\/}: #1\endaddress}
\fi

\topmatter
\title $q$-Rook polynomials and matrices over finite fields 
\endtitle
\rightheadtext{
$q$-Rook polynomials and matrices over finite fields}
\author James Haglund\endauthor
\date June 17, 1997 \enddate

\address Department of Mathematics, University of Illinois,
Urbana, Illinois 61801 \endaddress
\email jhaglund\@math.uiuc.edu \endemail

\thanks 
The author is supported by NSF grant DMS-9627432.  Part of this work was
done while the author was visiting MSRI during the Spring of 1997.
\endthanks
%  Math Subject Classifications 

\abstract
Connections between $q$-rook polynomials and matrices over finite
fields are exploited to derive a new statistic for Garsia and
Remmel's $q$-hit polynomial.  Both this new statistic $mat$ and 
another statistic for the $q$-hit polynomial $\xi$ recently
introduced by Dworkin are shown to induce different multiset Mahonian permutation
statistics for any Ferrers board.  In addition, for the triangular boards
they are shown to generate different families of Euler-Mahonian statistics.
For these boards the $\xi$ family includes Denert's statistic $den$, and gives
a new proof of Foata and Zeilberger's Theorem that 
$(exc,den)$ is jointly distributed with $(des,maj)$. 
The $mat$ family appears to be new.
A proof is also given that the $q$-hit polynomials are
symmetric and unimodal.
\endabstract
\endtopmatter

\document
\noindent
\head 1. Introduction \endhead
Notation:  LHS and RHS are abbreviations for ``left-hand-side" and
``right-hand-side", respectively.  $\Bbb N$ denotes the nonnegative
integers, $\Bbb Z$ the integers, $\Bbb P$ the positive integers,
and $\Bbb F_q$ a finite field with
$q$ elements.

A $board$ is a subset of an $n \times n$ grid of squares.  We
label the squares of the grid with the same (row,column) coordinates as the
squares of an $n \times n$ matrix; the
lower-left-hand-corner square has label $(n,1)$, etc. 
A $Ferrers$ $board$ is a board
with the property that 
$(i,j) \in B
\implies (k,p) \in B$ for $1\le k\le i$ and $j\le p\le n$.  Garsia
and Remmel [GaRe] introduced the following $q$-rook polynomial;

$$
R_k(B):=\sum_{C} q^{\text{inv}(C,B)}, \eqno(1)
$$
where the sum is over all placements $C$ of $k$ non-attacking rooks on the
squares of the Ferrers board $B$.  Non-attacking means no two rooks are 
in the same column, and no two are in the same row.  To calculate the
statistic $\text{inv}(C,B)$, cross out all squares which either contain a rook, or
are above or to the right of any rook.  The number of squares of $B$
not crossed out is $\text{inv}(C,B)$ (see Fig. 1).  

\input epsf
\midinsert
%make all figures 1/2 size
$$\vbox{\epsffile{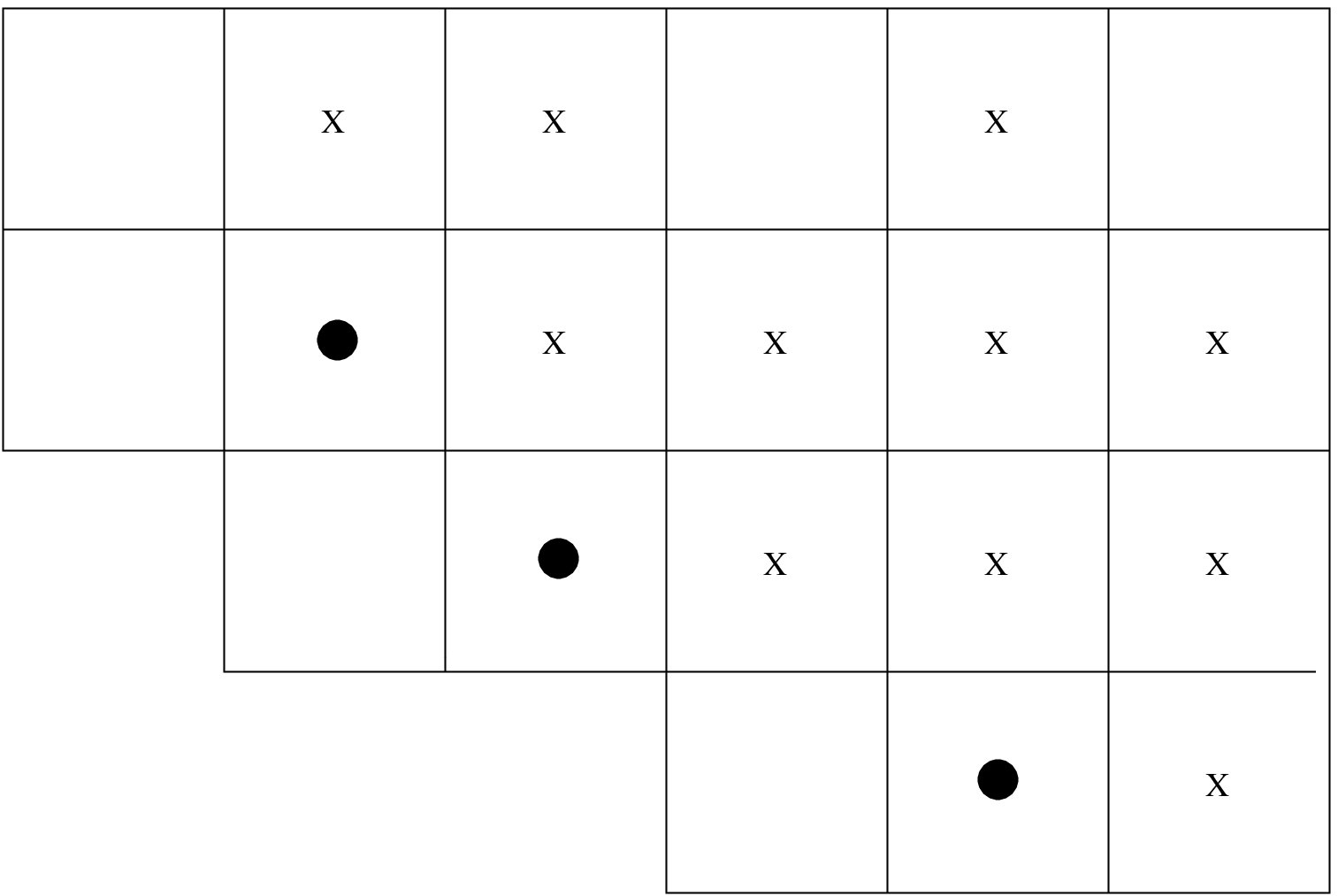}}$$
\botcaption{Figure 1}  A placement of $3$ rooks with $inv$ 
statistic $6$.
\endcaption
\endinsert

Garsia and Remmel showed that the $R_k$ enjoy many of the same
properties as the famous rook numbers $r_k$ introduced by Riordan
and Kaplansky [KaRi],[Rio].  For example,

$$
\sum_{k=0}^n [x][x-1]\cdots [x-k+1] R_{n-k}(B)
= \prod_{i=1}^n [x+c_i-i+1], \eqno(2)
$$
where $[x]:= (1-q^x)/(1-q)$ and $c_i:=$ the number of squares in the
$i$th column of $B$.  Note that our definition of a board requires
$c_n \le n$ (such boards are called $admissible$ in the literature).  This
assumption holds throughout the article, except as noted in Theorem 7.
 When $q \to 1$ in (2) we get a classic result of
Goldman, Joichi, and White [GJW].  As noted by Garsia and Remmel, an interesting
consequence of (2) is that two Ferrers boards have the same rook numbers if
and only if they have the same $q$-rook numbers, since both of these are
determined by the multiset whose elements are the shifted column
heights $c_i(B)-i+1$.

Letting $[k]!:=\prod _{1\le i\le k} [i]$,
and defining $T_k(B)$ via
$$
\sum_{k=0}^n [k]! R_{n-k}(B) \prod_{i=k+1}^n (x^k-q^i)
= \sum_{k=0}^n T_kx^k, \eqno(3)
$$
another result of Garsia and Remmel is that 
$$
T_k(B)=\sum_{C \atop \text{$n$ rooks, $k$ on $B$}}
q^{\text{stat$(C,B)$}},
$$
for some statistic $\text{stat$(C,B)$} \in \Bbb N$.  
In the sum above $C$ is a
placement of $n$ non-attacking rooks on the $n \times n$ grid,
with exactly $k$ on $B$.  For $q=1$ it reduces to $t_k(B)$, the
hit number of Riordan and Kaplansky, which equals the number of
permutations which ``hit" $k$ of the ``forbidden positions"
represented by the squares of $B$.  

Garsia and Remmel gave a recursive definition of $\text{stat$(C,B)$}$
, and left it as an open problem to determine a
method of generating $T_k(B)$ directly from the rook
placements (as in the definition of $R_k(B)$).  This problem
has recently been solved by M. Dworkin [Dwo], who shows that
$$
T_k(B)=\sum_{C \atop \text{$n$ rooks, $k$ on $B$}}
q^{\xi (C,B)},
$$
where $\xi(C,B)$ is calculated by the following procedure.

First place a bullet under each rook, and an $x$ to the right 
of any rook.  Next, for
each rook on $B$, place a circle in the empty cells of $B$ that are below 
it in the column.  Then for each rook off $B$, place a circle in the empty
cells below it in the column, and also in the empty cells of 
$B$ above it in the column.  Then $\xi (C,B)$ is the number of
circles.  See Fig. 2.

\input epsf
\midinsert
%make all figures 1/2 size
$$\vbox{\epsffile{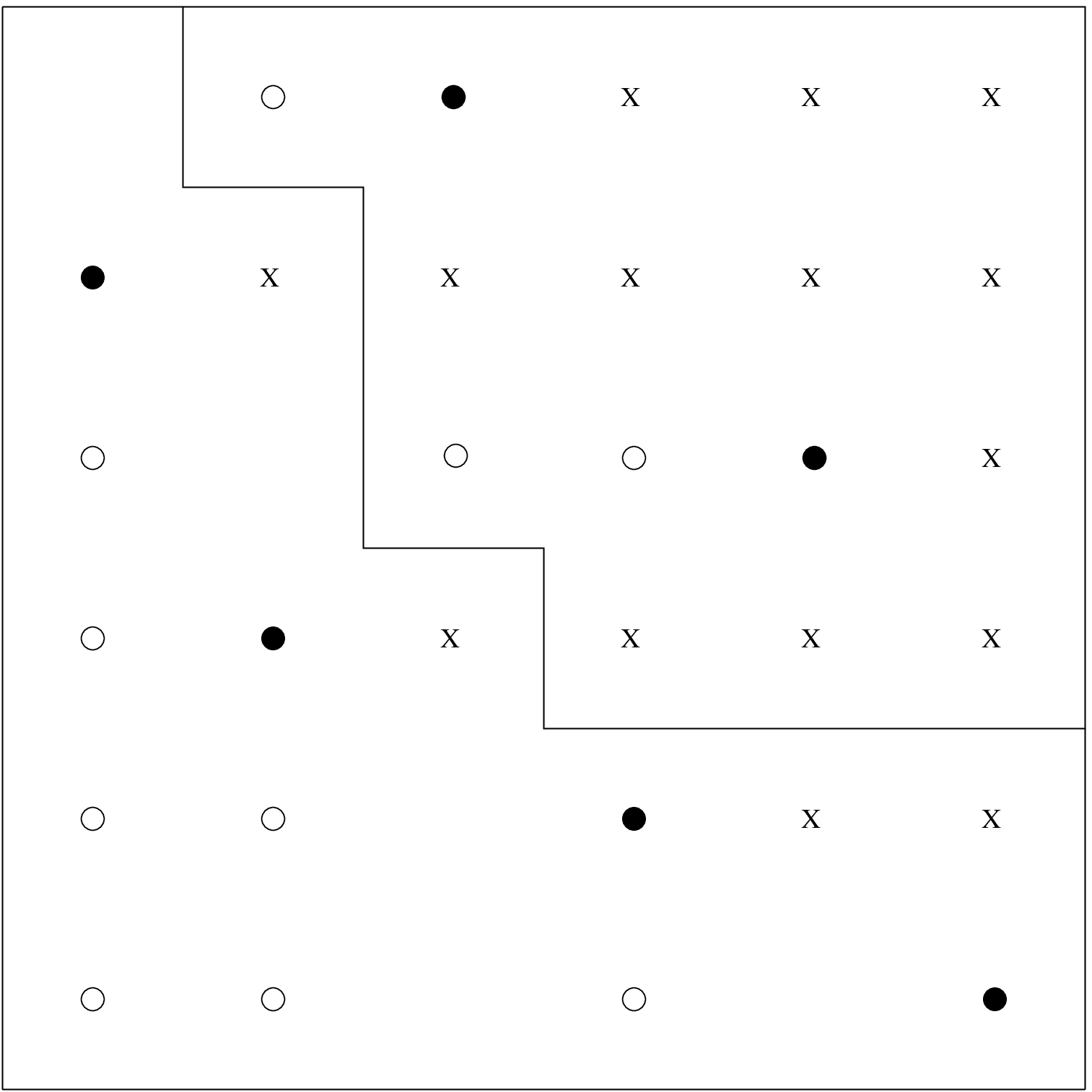}}$$
\botcaption{Figure 2}  A placement of $6$ rooks with $2$
rooks on $B$: $\xi=10$. 
\endcaption
\endinsert

The main result in this article is another solution to Garsia and Remmel's
problem, discovered before the author knew of Dworkin's result.
This new statistic, which
we call $mat$, bears
superficial similarities to Dworkin's $\xi$, but the author has been
unable to show that one being a solution implies the other is as
well. 
We arrive at $mat$ by counting matrices over finite fields
subject to certain constraints, while Dworkin first generalizes a
recurrence for the hit numbers given by Riordan, then
shows $\xi$ satisfies this recurrence.

A permutation $\sigma$ of a multiset $M$ 
is a linear list $\sigma_1 \sigma_2 \cdots \sigma _{\#M}$ 
of the elements of $M$.
For any vector $\bold v=(v_1,v_2,\ldots ,v_t)$ of nonnegative integers,
let $\{1^{v_1}2^{v_2}\cdots t^{v_t}\}$ denote the multiset having
$v_i$ copies of $i$, and let $M(\bold v)$ be the set of 
permutations of 
$\{1^{v_1}2^{v_2}\cdots t^{v_t}\}$. If $\bold v =(1,1,\ldots ,1)$ is 
the vector with $n$ ones, we identify the element
$\sigma_1 \sigma_2 \cdots \sigma_n \in M(\bold v)$ with the
element $\left (\matrix 1 & 2& \cdots &n\\
\sigma_1 & \sigma_2 & \cdots &\sigma_n \endmatrix \right )$ of the symmetric
group $S_n$.

A statistic $stat$ on permutations in $S_n$
is called $Mahonian$ if
$$
\sum_{\sigma \in S_n} q^{\text{stat}(\sigma)}=[n]!.
$$
It is called {\it multiset Mahonian} if
$$
\sum_{\sigma \in M(\bold v)} q^{\text{stat}(\sigma)}=
\left[ \matrix & \#M & \\ v_1,& v_2, & \ldots &,v_t \endmatrix \right ]
$$
for all vectors $\bold v$, where 
$\left[ \matrix & \#M \\ v_1,& v_2, & \ldots &,v_t \endmatrix \right ]
:=\frac{[\#M]!}{\prod_{i=1}^t[v_i]!}$ is the $q$-multinomial coefficient.
The study of Mahonian statistics has become a large enterprise in 
recent years.  
Dworkin showed that $\xi$ induces a Mahonian statistic for any Ferrers
board $B$, and we generalize this to show how $\xi$ and $mat$ both
induce multiset Mahonian statistics.  We should mention that Dworkin
gave his definition and results for $\xi$ in the more general setting of
$skyline$ $boards$, which
are obtained by permuting the columns of a Ferrers board.
  Unfortunately,
if we extend our definition of $mat$ in a straightforward way, the
resulting statistic is not Mahonian for skyline boards; in fact, it
is not even nonnegative.  For that reason, we will restrict our 
attention to Ferrers boards in this article.

For $\sigma \in M(\bold v)$,  
a descent of $\sigma$ is a
value of $i$, $1\le i< n$, such that $\sigma_i > \sigma_{i+1}$, where
$n=\sum_iv_i$.
MacMahon showed the statistic $maj$ is multiset Mahonian, where
$$
\text{maj}(\sigma):=\sum_{i:\, \sigma_i > \sigma_{i+1}}i.
$$
Let $des(\sigma)$ denote the number of descents of $\sigma$.
A pair $(stat1,stat2)$ of statistics on permutations in $S_n$ 
is called $Euler$-$Mahonian$ if it is
jointly distributed with $(des,maj)$, i.e. if
$$
\sum_{\sigma \in S_n } 
p^{\text{stat1}(\sigma)}
q^{\text{stat2}(\sigma)}=
\sum_{\sigma \in S_n}
p^{\text{des}(\sigma)}
q^{\text{maj}(\sigma)}.
$$
Dworkin noted that $\xi(B)$ is part of an Euler-Mahonian 
pair if $B$ is a triangular-shaped board.
In section $4$ we develop this idea further, and show how 
$\xi$ and $mat$ both induce families of eight
Euler-Mahonian pairs of statistics each of the form $(des,stat)$, which
are all different from
one another. 

Let $f(\bold v):=f_1(\bold v)\cdots f_n(\bold v)$ be the unique element 
of $M(\bold v)$ with no descents.  An $excedence$ of 
$\sigma \in M(\bold v)$ is a value of $i$ such that
$\sigma_i > f_i(\bold v)$, and we denote the number of such
excedences by $\text{exc}(\sigma)$.  For example, the permutation
$\sigma = 2313212$ has $3$ excedences, occurring in the first,
second, and fourth places of $\sigma$, and so 
$\text{exc}(\sigma)=3$.

Not many Euler-Mahonian pairs are known.
A general overview of the few that have been discovered
can be found in [CSZ].  There the authors 
classify a solution to the equation
$$
\sum_{\sigma \in S_n \atop \text{exc}(\sigma)=k}
q^{\text{statx} (\sigma)} =
\sum_{\sigma \in S_n \atop \text{des}(\sigma)=k}
q^{\text{maj} (\sigma)}, \eqno(4)
$$
as a ``proper" Euler-Mahonian pair, proper
indicating that $exc \ne des$.
Using the geometry of the board, it is simple to convert 
one of our Euler-Mahonian pairs $(des,stat)$ into a solution to (4).
When this conversion is applied to the $\xi$ family, we get  
Denert's statistic $den$ [Den], and
a new proof of a result of Foata and Zeilberger [FoZe],
that $(exc,den)$ is jointly distributed with $(des,maj)$.
On the other hand, the $mat$ family gives what appears to be a
fundamentally new solution to (4).

Garsia and Remmel also obtained a solution to a form of (4), namely
$$
\sum_{\sigma \in S_n \atop \# \{i:\, \sigma_i \ge i \}=k}
q^{\text{naj} (\sigma)} =
\sum_{\sigma \in S_n \atop \text{des}(\sigma)=k}
q^{\text{maj} (\sigma)},
$$
involving a statistic $naj$ which they defined recursively.  They also
gave a recursive definition of a Mahonian statistic which involved an
arbitrary Ferrers board.  It would be interesting to obtain 
non-recursive versions of the definitions of these statistics, and
determine how they relate to other Mahonian statistics and Euler-Mahonian
pairs.

In section $5$ we show that $T_k(B)$ is a
symmetric and unimodal polynomial in $q$ for all $B$, a fact first
proved in [Ha1].  The proof is a simple extension of
Garsia and Remmel's proof that $T_k(B) \in \Bbb N[q]$.  
For some boards we prove a stronger result by
a different method.

\head 2. Matrices over Finite Fields\endhead
\noindent

Solomon [Sol] showed how a placement of $k$ non-attacking
rooks on a rectangular board can naturally be
associated to a rectangular $n \times m$
matrix with entries in $\Bbb F_q$ and of rank $k$.
Ding has shown that a similar construction involving
matrices over the complex numbers in the shape of a
Ferrers board has applications to
topological questions involving certain algebraic varieties
[Din1],[Din2].
In the lemma below we generalize Solomon's result to Ferrers boards;
the proof is a straightforward extension of his.

\proclaim{Definition 1}  For $B$ a Ferrers board with $n$
columns (some of which may be empty), let $P_k(B)$ be the number of
$n \times n$ matrices $A$ with entries in $\Bbb F_q$, of rank $k$, and 
with the restriction that all the entries of $A$ in those squares of $A$
outside of $B$ are zero.  For example, if $B$ is the board
consisting of squares $(1,2)$,$(1,3)$, and $(2,3)$, then
$P_0=1$,$P_1=2q^2-q-1$,$P_2=q(q-1)^2$, and $P_3=0$.
\endproclaim

\proclaim{Theorem 1}  For any Ferrers board $B$,
$$
P_k(B)=(q-1)^k q^{\text{Area}(B)-k}R_k(q^{-1}),
$$
where $Area(B)$ is the number of squares of $B$.
\endproclaim
\noindent 
$Proof:$ Let $A$ be a matrix of rank $k$, with entries in $\Bbb F_q$,
and zero outside of $B$.  We perform an operation on $A$
which we call the $elimination$ $procedure$.  Starting at the
bottom of column 1 of $A$, travel up until you arrive at a
nonzero square $\beta$ (if the whole first column is zero go to column 2
and iterate).  Call this nonzero square a {\it pivot spot}.  Next
add multiples of the column containing $\beta$ to the columns to the
right of it to produce zeros in the row containing $\beta$ to the
right of $\beta$.  Also add multiples of the row containing $\beta$
to the rows above it to produce zeros in the column containing
$\beta$ above $\beta$.  Now go to the bottom of the next column and
iterate; find the lowest nonzero square, call it a pivot spot,
then zero-out entries
above and to the right as before.

If we place rooks on the square $\beta$ and the other pivot spots
we end up with $k$ non-attacking rooks.  The number of matrices which
generate a specific rook placement $C$ is 
$$
(q-1)^k q^{\text{\# of squares to the right of or above a rook}}
$$
$$
= (q-1)^k q^{\text{Area}(B)-k-\text{inv}(C,B)}. \qquad \blacksquare
$$

\proclaim{Corollary 1}  Let $P_k$ be the number of $n \times n$ upper
triangular matrices of rank $k$ with entries in $\Bbb F_q$.
Then 
$$
P_k = (q-1)^k q^{{n+1 \choose 2}-k}S_{n+1,n+1-k}(q^{-1}),
$$
where $S_{n,k}(q)$ is the $q$-Stirling number of the second kind defined by
the recurrences
$$
S_{n+1,k}(q):=q^{k-1}S_{n,k-1}(q) + [k]S_{n,k}(q) \qquad (0\le k\le n+1),
$$
with the initial conditions $S_{0,0}(q)=1$ and $S_{n,k}(q)=0$ for
$k<0$ or $k>n$.
\endproclaim
\noindent
$Proof:$ It is known [GaRe,p.248] that if $B$ is the triangular board whose
$i$th column has height $i$, then
$$
R_k(B)=S_{n+1,n+1-k}(q).
$$
Now apply Theorem 1. $\qquad \blacksquare$

\proclaim{Corollary 2}  For any Ferrers board $B$,
$$
\sum_{k=0}^n (1-x)(1-xq)\cdots (1-xq^{k-1}) P_{n-k}(B)
=\prod_{i=1}^n(q^{c_i}-xq^{i-1}).
$$
\endproclaim
\noindent
$Proof:$  This is obtained by replacing by replacing $q$ by $q^{-1}$ in (2),
applying Theorem 1, and doing other simple transformations such as
replacing $q^x$ by $1/x$. $\qquad \blacksquare$

\noindent $Remark:$  In [Hag1], the following identity was derived as a
limiting case of a hypergeometric result:
$$
\sum_k R_k(B)(1-q)^k =1
$$
(this can also be obtained by letting $x \to \infty$ in (2)).
Using Theorem 1, this is equivalent to the trivial statement
$$
\sum_k P_k(B)=q^{\text{Area}(B)}.
$$

\proclaim{Definition 2}  Let $C$ be a placement of $n$ non-attacking
rooks on the $n \times n$ grid, with $k$ rooks on the Ferrers board $B$.
Define $\text{cross}(C,B)$ 
to be the number of squares of the $n \times n$
grid satisfying one of the following conditions:

- containing a rook or to the right of   
a rook

- above a rook and on $B$

- below a rook which is off $B$

\noindent Furthermore let $\text{mat}(C,B):=n(n-k) +\text{Area}(B)-
\text{cross}(C,B)$.
  See Fig. 3.
\endproclaim

\input epsf
\midinsert
%make all figures 1/2 size
$$\vbox{\epsffile{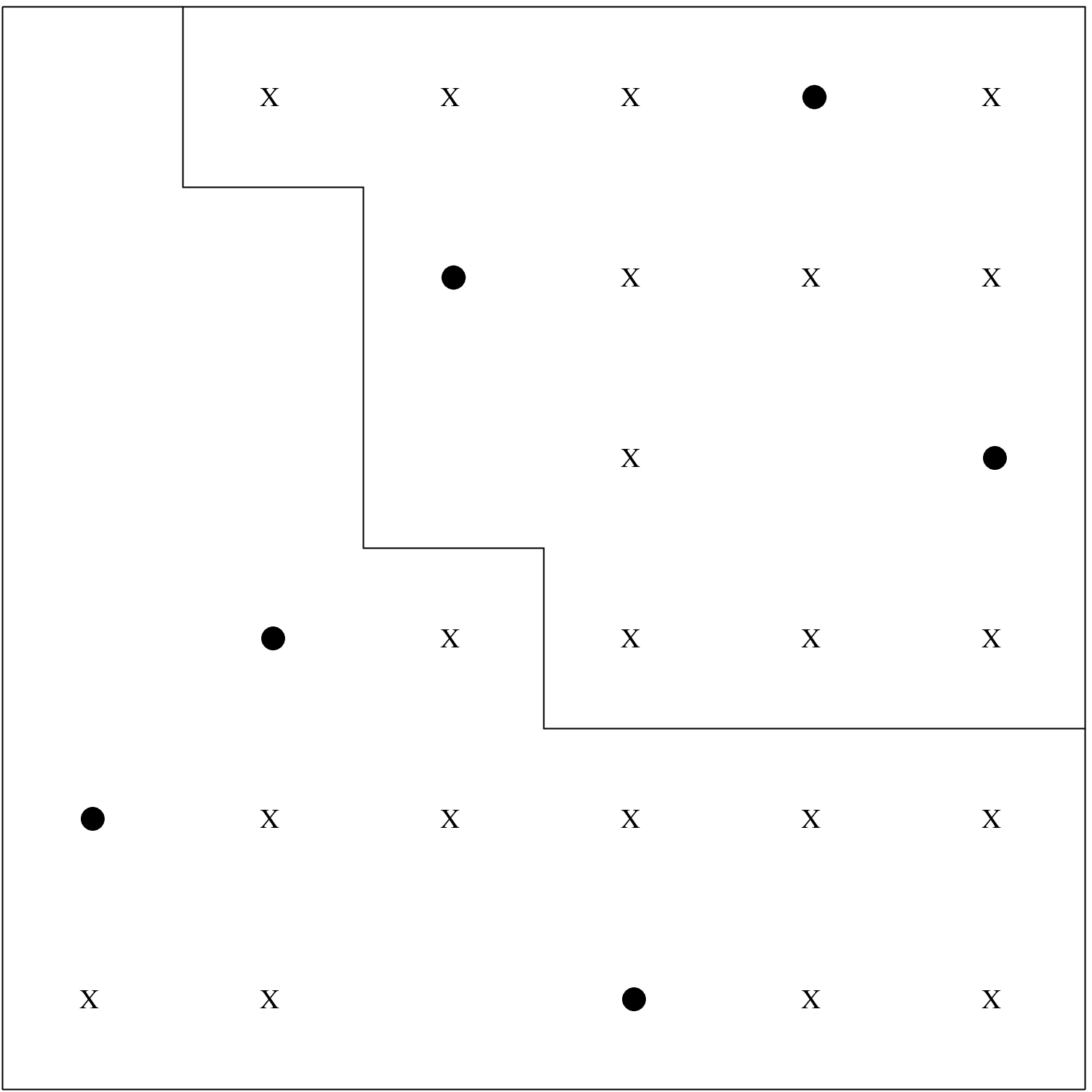}}$$
\botcaption{Figure 3}  A placement of $6$ rooks with $3$
rooks on $B$: $cross=27$ so 
${mat}=3*6+16-27=7$. 
\endcaption
\endinsert

\proclaim{Theorem 2}  If $B$ is any Ferrers board,
$$
T_k(B)=\sum_{C \atop \text{$n$ rooks, $k$ on $B$}} q^{\text{mat}(C,B)}.
$$
\endproclaim
\noindent
$Proof:$  Replacing $q$ by $q^{-1}$ in (3) and
multiplying by $(q-1)^n q^{\text{Area}(B)}$ we get
$$
\sum_{k=0}^n (q^k-1)\cdots (q^2-1)(q-1)(q-1)^{n-k}q^{\text{Area}(B)-
(1+2+\ldots +k-1 +k+1+\ldots +n)} R_{n-k}(q^{-1})
$$
$$
\times \prod_{i=k+1}^n (q^ix-1)= \sum_{k=0}^n x^k (q-1)^n
q^{\text{Area}(B)}T_k(q^{-1}),
$$
or
$$
\sum_{k=0}^n (q^k-1)\cdots (q-1)(q-1)^{n-k}q^{\text{Area}(B)}
R_{n-k}(q^{-1})q^{-(n-k)}
\prod_{i=k+1}^n (q^ix-1)
$$
$$
= \sum_{k=0}^n x^k (q-1)^n
q^{\text{Area}(B)+{n \choose 2}}T_k(q^{-1}), \eqno (5)
$$
or
$$
\sum_{k=0}^n (q-1)(q^2-1)\cdots (q^k-1)P_{n-k}\prod_{i=k+1}^n
(q^ix-1) = \text{ RHS of (5) }=\sum_{k=0}^n x^k Q_k \eqno(6)
$$
say.  We will prove Theorem 2 by showing that
$$
Q_k = \sum_{C \atop \text{$n$ rooks, $k$ on $B$}}
(q-1)^n q^{\text{Area}(B) + {n \choose 2}-\text{mat}(C,B)}
$$
$$
= \sum_{C \atop \text{$n$ rooks, $k$ on $B$}}
(q-1)^n q^{{n \choose 2}+\text{cross}(C,B)-n(n-k)}. \eqno(7)
$$
Our strategy will be to exploit the combinatorial interpretation of the
LHS of (6).  Using the following special case of Cauchy's famous
$q$-$binomial$ $theorem$;
$$
\prod_{j=0}^{m-1}(1+xq^j) = \sum_{k=0}^m 
\left[ \matrix m \\ k \endmatrix \right ]
q^{{k \choose 2}}x^k, \eqno(8)
$$
where 
$\left[ \matrix m \\ k \endmatrix \right ]:=\frac{[m]!}{[k]![m-k]!}$
is the $q$-binomial coefficient, the coefficient of $x^s$ in the LHS of
(6) can be written as
$$
\sum_{k=0}^{n-s}P_{n-k}(q^n-q^{n-k})(q^{n-1}-q^{n-k})\cdots
(q^{n-k+1}-q^{n-k})
$$
$$
\times
\left[ \matrix n-k \\ s \endmatrix \right ](-1)^{n-k-s}
q^{{s \choose 2} +s(k+1)-k(n-k)}.  \eqno(9)
$$
We want to show that the expression above equals the RHS of (7).  Let $A$ be
a matrix of rank $n-k$, with entries in $\Bbb F_q$ and zero outside $B$.
We now 
perform an operation on $A$ which we call the {\it 
replacement procedure}.  Starting with the last row (the bottom row) of $A$,
define row $\alpha _k$ as the bottom-most row linearly dependent on the rows
below it (or $\alpha _k=n$ if the last row is zero).  Next let 
row $\alpha _{k-1}$ be the next bottom-most row 
linearly dependent on the rows below it, etc.  Thus
we end up with $k$ rows $\alpha _k > \alpha _{k-1} > \cdots > \alpha _1$.
We call the rows $\alpha _1,\ldots ,\alpha _k$ ``dependent rows" and the
other rows of $A$ ``keeper rows".

Now replace row $\alpha _1$ by any of the $q^n-q^{n-k}$ rows which are
linearly independent of the rows of $A$.  
Call this new row $\overline{\alpha _1}$,
and note that $P_{n-k}$ is multiplied by $(q^n-q^{n-k})$ in 
(9).  
If $\overline{\alpha _1}$ has any nonzero entries off $B$, we call it a 
{\it pivot row}, and the spot where the left-most nonzero entry in
$\overline{\alpha _1}$ occurs a {\it pivot spot} ( 
when we perform the elimination
procedure later, this spot will be a pivot).  Next replace row 
$\alpha _2$ by a new row linearly independent of both the 
rows of $A$ and
the new row $\overline{\alpha _1}$, 
with the added constraint that if $\overline{\alpha _1}$ is a 
pivot row, we require the new $\overline{\alpha _2}$ row to have a zero in the column
containing the pivot spot in row $\overline{\alpha _1}$.  If $\overline{\alpha _1}$ is a pivot row,
there are $q^{n-1}-q^{n-k}$ choices for $\overline{\alpha _2}$ (if we look at all
linear combinations of $\overline{\alpha _1},w_1,\ldots ,w_{n-k}$
where the $w_i$ are the keeper rows of $A$, then for any fixed
$c_2,\ldots ,c_{n-k+1}$, the sums
$$
c_1\overline{\alpha _1}+c_2w_1+\ldots +c_{n-k+1}w_{n-k}
$$
produce $q$ different values in the column containing the pivot spot
of $\overline{\alpha _1}$ as $c_1$ cycles through its $q$ possible values) and
$q^n-q^{n-k+1}$ choices otherwise.  In the latter case, we define the
{\it weight} of row $\overline{\alpha _2}$ to be $q^{-1}$ (this is what we need to
multiply $q^n-q^{n-k+1}$ by to get the desired factor $q^{n-1}-q^{n-k}$
occurring in (9)).  If $\overline{\alpha _1}$ is a pivot row, let the weight of
row $\overline{\alpha _2}$ be 1.  
As before, if $\overline{\alpha _2}$ contains any nonzero entries off
$B$ we call it a pivot row, and its left-most nonzero entry a pivot spot.

Now for $\overline{\alpha _3}$, we require there be zeros in the columns containing any
pivot spots in rows $\overline{\alpha _1}$ or $\overline{\alpha _2}$.  More generally, in $\overline{\alpha _j}$,
we require zeros below any of the pivot spots in rows $\overline{\alpha _i}$, 
$1\le i <j$, and define the weight of $\overline{\alpha _j}$ to be
$q^{-w}$, with $w$ equal to the cardinality of $\{ {\overline{\alpha _i}}:\, 
{\overline{\alpha _i}} \text{ is not a pivot row and } 1\le i <j\}$.  

Let 
$$
\text{factor}(s,k):=
\left[ \matrix n-k \\ s \endmatrix \right ](-1)^{n-k-s}
q^{{s \choose 2} +s(k+1)-k(n-k)}
$$
(as in (9)).  The argument above shows that $Q_s$ equals the
number of matrices (counted with multiplicity) which are
obtained by starting with matrices which are zero
outside of $B$ and
performing the replacement procedure, and finally multiplying by the appropriate
weight and factor.  If we perform the elimination procedure from the proof
of Theorem 1 to one of these new matrices, we end up with $n$ pivots, where
the pivots off $B$ are exactly those pivots spots defined above from the
$\overline{\alpha _i}$.

Let $Q_{s,j}$ be the number of these matrices, counted with weights and
factors, with $j$ pivots off $B$.

\vskip .1in
\noindent 
Case 1: $j=n-s$.  In this case, all the weights are 1.  The row numbers
with pivots off $B$ must have been the original $\alpha _1 < \alpha _2 <
\cdots < \alpha _{n-s}$.  There are $q^{n-\alpha _{n-s}}$ choices of row
$\alpha _{n-s}$ to be dependent on the rows below; we call 
$q^{n-\alpha _{n-s}}$ the {\it pre-image term} for this row.  There are
$q^{n-\alpha _{n-s-1}-1}$ (the pre-image term for this row) choices for 
row $\alpha _{n-s-1}$ to be 
dependent on the rows below, etc.  Note that the $P_i$ in (9)
satisfy $i\ge s$, and only the $i=s$ term can possibly generate matrices
with $n-s$ pivots off $B$.  By the elimination procedure and (9),
$$
Q_{s,n-s} = \sum_{C \atop \text{$n$ rooks, $s$ on $B$}}
(q-1)^n q^{n-\alpha _{n-s}+n-\alpha _{n-s-1} -1 +\ldots +
n-\alpha _{1} -(n-s-1)}
$$
$$
\times q^{\text{\# of squares to the right of a rook, or above a rook
and on $B$}}
$$
$$
\times q^{-\text{\# of squares below a rook off $B$ 
and to the right of some rook}}
$$
$$
%\times q^{{s \choose 2} +s(n-s+1)-(n-s)s},
\times \text{factor}(s,n-s),
$$
where in the sum above $\alpha _1(C) < \cdots < \alpha _{n-s}(C)$ are the
row numbers with rooks off $B$.
Now 
$$
\multline
\text{\# of squares to the right of a rook or above a rook and on $B$}\\ -
\text{\# of squares below a rook off $B$ and to the right of some
rook}\\ +
n-\alpha_{n-s}+n-\alpha_{n-s-1}+\ldots +n-\alpha_1 =
\text{cross}(C,B)-n,
\endmultline
$$
and plugging this in above, after a short calculation we get
$$
Q_{s,n-s} = \sum_{C \atop \text{$n$ rooks, $s$ on $B$}}
(q-1)^n q^{{n \choose 2}+\text{cross}(C,B)-n(n-s)}.
$$
In view of (7), Theorem 2 follows if we can show $Q_{s,j}=0$ if
$j<n-s$.

\vskip .1in
\noindent
Case 2: $j<n-s$. By an abuse of terminology, if the weight of a
row is $q^{w}$, we sometimes refer to $w$ as the weight.
A similar remark applies to the pre-image term.

For each $k$ with $j\le k\le n-s$, the term $P_{n-k}$
in (9) makes a contribution.  Say after replacement
and elimination, we end up with a placement $C$ of $n$ rooks,
with the $j$ pivots off $B$ in rows
$\beta _1 < \beta _2 < \cdots < \beta _j$.  Then all these rows, and 
$k-j$ others, must have been the original $\alpha _1, \ldots ,
 \alpha _k$.  We have to sum over all choices of the $k-j$ others,
 taking into account the weights, the pre-image terms
 $n-\alpha _k + n-\alpha _{k-1}-1 +\ldots +n-\alpha _1-(k-1)$, and the
$\text{factor}(s,k)$ term from (9).  

Say there are
$\mu _0$ new rows above row $\beta _1$, $\mu _1$ new rows
between rows $\beta _1$ and $\beta _2$,$\ldots$, and $\mu _j$ below  
row $\beta _j$, with $\mu _i \ge 0$ and 
$\mu _0 +\mu _1 +\ldots + \mu _j=k-j$.  Lets compute the
total weight of such an arrangement, using the fact that as we
move downwards, the weights of the rows decrease by one each time,
unless the row is just below a $\beta _i$, in which case the weight stays
the same.  

$$ \text{The $\mu _0$ rows above row $\beta _1$ have weights
$-0,-1,\ldots ,-\mu _0+1$}.
$$
$$ \text{Row $\beta _1$ has weight
$-\mu _0$}.
$$
$$ \text{The $\mu _1$ rows between rows $\beta _1$ and $\beta _2$
have weights
$-\mu _0, -\mu _0-1,\ldots ,-\mu _0-\mu _1+1$}.
$$
$$ \text{Row $\beta _2$ has weight
$-\mu _0-\mu _1$}.
$$
$$
\vdots
$$
$$ \text{Row $\beta _j$ has weight
$-\mu _0-\mu _1-\ldots -\mu _{j-1}$}.
$$
$$ 
\medmuskip0mu
\text{The $\mu _j$ rows below row $\beta _j$ have weights
$-\mu _0 - \mu _1 -\ldots -\mu _{j-1}, \ldots ,
-\mu _0 - \mu _1 -\ldots -\mu _{j}+1$}.
$$

For the $\mu _j$ rows below row $\beta _j$, say rows
$\gamma _1 > \gamma _2 > \cdots > \gamma _{\mu _j}$, the
sum of the pre-image terms will be
$n-\gamma _1 + n-\gamma _2 -1 +\ldots + n-\gamma _{\mu _j}
- \mu _j+1$, and combining this with the weight for these
rows gives a total contribution of
$$
\multline
n-\gamma _{\mu _j}-\mu _j +1-\mu _0 -\mu _1 -\ldots -\mu_{j-1}
+ n-\gamma _{\mu_{j}-1}-(\mu _j-2)-\mu _0 -\mu _1 - \\ \ldots
-\mu_{j-1}-1 +\ldots +n-\gamma _1-\mu _0-\mu _1-\ldots -\mu _j+1,
\endmultline
$$
or
$$
n-\gamma _1 +n-\gamma _2+\ldots +n-\gamma _{\mu _{j}}
-\mu _j(k-j-1).
$$

Lets skip past row $\beta _j$ for the moment and consider rows
$$
\gamma _{\mu_j+1} > \gamma _{\mu _j +2} > \cdots >
\gamma _{\mu _j+\mu _{j-1}}
$$ between rows $\beta _j$ and $\beta _{j-1}$.  The pre-image terms
will be
$$
n-\gamma _{\mu_j+1} -(\mu _j+1) + n-\gamma _{\mu _j +2} -(\mu _j +2)
+ \ldots +
n-\gamma _{\mu _j+\mu _{j-1}} - (\mu _j + \mu _{j-1}).
$$
Adding in the weights as before we end up with a contribution of
$$
n-\gamma _{\mu_j+1} -1 + n-\gamma _{\mu _j +2} -1
+ \ldots +
n-\gamma _{\mu _j+\mu _{j-1}} - 1 - \mu _{j-1}(k-j-1)
$$ for these rows.  
Continuing in this way,
for the $\mu _0$ rows above row $\beta _1$ we get a total
contribution of
$$
n-\gamma _{\mu_j+\ldots +\mu _1+1} -j  
+ \ldots +
\gamma _{\mu _j+\ldots +\mu _{0}} - j - \mu _{0}(k-j-1).
$$

As we range over all legal choices of the $\gamma _i$ (i.e. $\gamma_i
\ne \beta_k$ for all $i,k$),
the numbers 
$$
\multline
n-\gamma_1,\ldots,n-\gamma_{\mu_j},n-\gamma_{\mu_j+1}-1,\ldots ,
n-\gamma_{\mu_j+\mu_{j-1}}-1,\\ \ldots ,
n-\gamma_{\mu_j+\ldots +\mu_1+1}-j,\ldots, 
n-\gamma_{\mu_j+\ldots +\mu_0}-j 
\endmultline
$$
range over all numbers between $0$ and $n-j-1$.  Thus raising $q$
to the power
of all terms above (ignoring the $\beta_i$ weights) gives 
$$
q^{-(k-j-1)(\mu _0+\ldots +\mu _j)} \times 
\text{the coefficient of $x^{k-j}$ in
$\prod _{i=0}^{n-j-1}(1+xq^i)$}. \eqno(10)
$$
By (8), (10)
reduces to 
$$
q^{-(k-j-1)(k-j)}
\left[ \matrix n-j \\ k-j \endmatrix \right ]
q^{{k-j \choose 2}}.
$$

Next we add in the contribution from the $\beta _i$.
The weight of $\beta _j$ is $-(\mu _0 +\ldots +\mu _{j-1})$
and its pre-image term is $n-\beta _j-\mu _j$.  
For $\beta _2$, the weight is $-(\mu _0 +\mu _1)$ and the 
pre-image term is $n-\beta _2 -(j-2)-(\mu _j +\ldots +\mu _2)$.
For $\beta _1$, the weight is $-\mu _0$ and the 
pre-image term is $n-\beta _1 -(j-1)-(\mu _j +\ldots +\mu _1)$.
The total contribution from the $\beta _i$ is thus
$q^{\beta}$, where
$$
\beta := n-\beta _1+\ldots +n-\beta _j -j(\mu _0 +\ldots +\mu _j)
-{j \choose 2}.
$$
For fixed $C$ (which also fixes the $\beta_i$) we thus have a contribution
to $Q_{s,j}$ of
$$
\multline
\sum_{k}(q-1)^n 
\left[ \matrix n-j \\ k-j \endmatrix \right ]
q^{{k-j \choose 2}-(k-j-1)(k-j)-j(k-j) -{j \choose 2}}
q^{n-\beta _1 +\ldots +n-\beta _j} \\
q^{\text{\# of squares to the right of a rook, or above a
rook and on $B$}} \\
q^{-\text{\# of squares below any rook off $B$ and to the right
of some rook}} \\
\times q^{{s \choose 2}+s(k+1)-k(n-k)}
\left[ \matrix n-k \\ s \endmatrix \right ]
(-1)^{n-k-s} \endmultline
$$
$$
=q^{\text{cross}(C,B) -n + {s+1 \choose 2} }
 \sum_k 
\left[ \matrix n-j \\ k-j \endmatrix \right ]
\left[ \matrix n-k \\ s \endmatrix \right ] (-1)^k
q^{k(k/2+s-n+1/2)}.
$$
Letting $k=n-s-u$, the sum above
reduces to 
$$
d\left[ \matrix n-j \\ s \endmatrix \right ]
\sum_{u\ge 0} 
\left[ \matrix n-j-s \\ u \endmatrix \right ]
(-1)^u
q^{{u \choose 2}}
$$
(with $d$ independent of $u$) which equals zero for $j<n-s$ by (8).
This completes the proof of Theorem 2. $\qquad \blacksquare$

\head 3.  Multiset Mahonian Statistics \endhead

A placement $C$ of $n$ rooks on an $n \times n$ grid can be identified with
a permutation $\sigma_1 \sigma_2 \cdots \sigma_n \in S_n$, 
in a simple way: 
a rook is on square $(i,j)$ if and only if  
$\sigma _i=j$.  We call this placement $C(\sigma)$ the $graph$ of $\sigma$.
Hence both $\xi$ and $mat$ can be regarded as
permutation statistics if we define $\xi(\sigma,B):=\xi(C(\sigma),B)$
and $\text{mat}(\sigma,B):=\text{mat}(C(\sigma),B)$.  

%If $B$ consists of the whole
%$n \times n$ grid, $mat(\sigma,B)$ equals the number of inversions of 
%$\sigma$.

Dworkin proved that $\xi$ is Mahonian for all
Ferrers boards $B$, i.e.
$$
\sum_{\sigma \in S_n} q^{\xi (\sigma,B)} = [n]!,
$$
or equivalently
$$
\sum_k T_k(B) = [n]!. \eqno (11)
$$
Eq. (11) can also be obtained by letting $x \to \infty$ in the
following known formula [Hag1,p.100], [Dwo,pp.35,38] 
$$
\sum_k \left [ \matrix x+k \\ n \endmatrix \right ] T_k(B) =
\prod_{i=1}^n [x+c_i-i+1],
$$
which in turn follows from (21) and another form of the $q$-binomial theorem.

In this section we show how to construct multiset Mahonian statistics from both
$\xi$ and $mat$.  
\proclaim {Definition 3} Let $B$ be a Ferrers board.  A $section$ $D$ of $B$ of
width $d$ is a set of $d$ consecutive columns of the $n \times n$ grid with the
property that the height of all the columns of $B$ in $D$ is the same.
%$c_i(B)$ is constant for all $i$ such that $\text{column} \,i \in D$.  
Say we have a placement $C$ of $n$ rooks on the
$n \times n$ grid, with $s$ of the $d$ rooks in $D$ on $B$.  We say $C$ is $D$-$standard$ 
if both of the following hold.
\vskip .1in
\noindent 
1) The $d-s$ rooks off $B$ and in $D$ are in the $d-s$ left-most columns of $D$.
   Furthermore, these $d-s$ rooks are in ``descending" order; if two of these
   rooks occupy squares $(i,j)$ and $(k,l)$, with $i< k$, then $j< l$.
\vskip .1in

\noindent
2) The $s$ remaining rooks in $D$ and on $B$ are in ``ascending" order; if
   two of these rooks occupy squares $(i,j)$ and $(k,l)$, with $i< k$,
   then $j > l$.
\vskip .1in
\noindent
Call a placement of $n$ rooks on the
$n \times n$ grid $B$-$standard$ if it is 
$D$-standard for all possible sections 
$D$ of $B$.  See Fig. 4.

\input epsf
\midinsert
%make all figures 1/2 size
$$\vbox{\epsffile{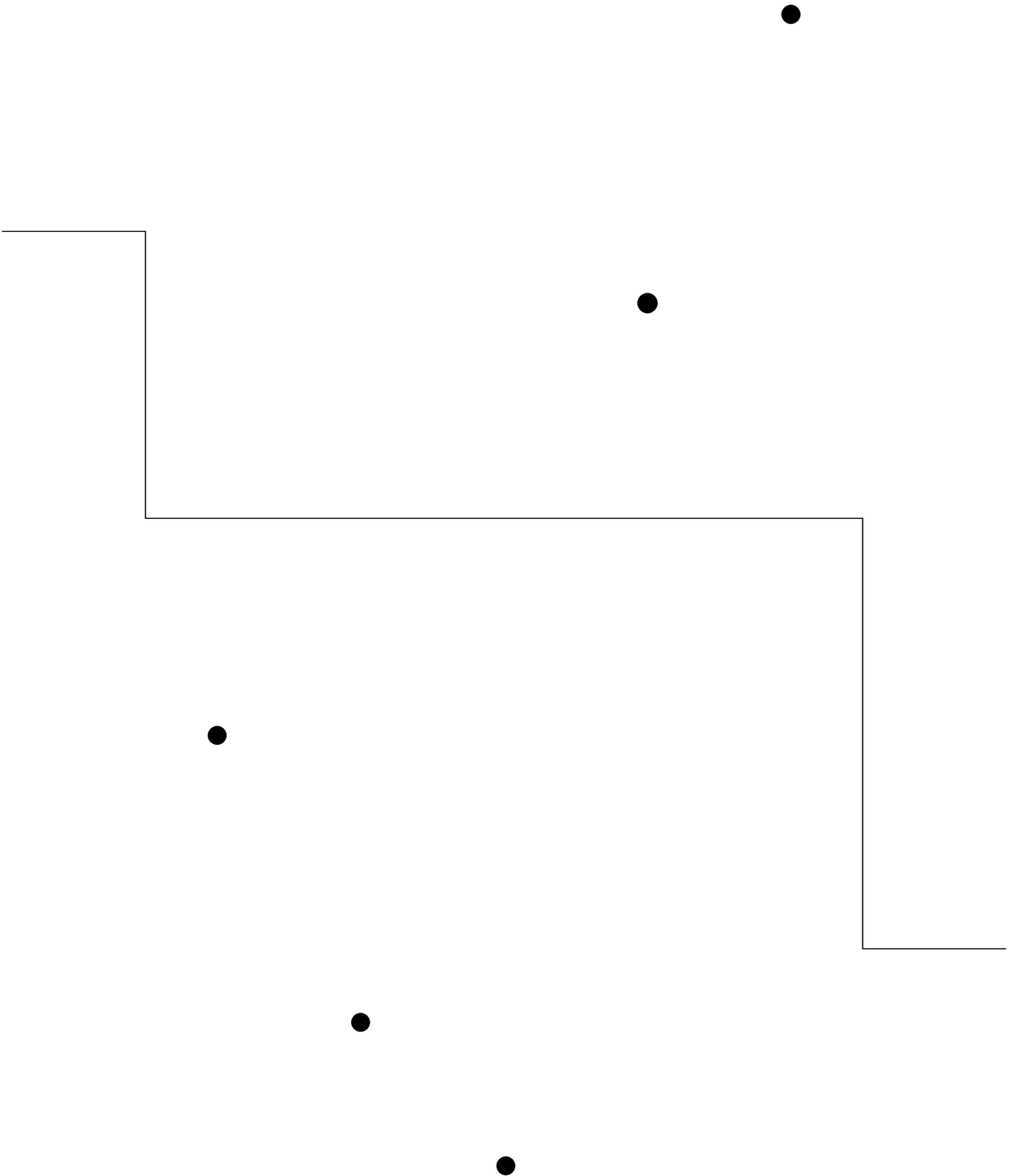}}$$
\botcaption{Figure 4}  A section $D$ of width $5$ and a 
$D$-standard placement of rooks (rooks outside of $D$ are
not pictured).
\endcaption
\endinsert

\endproclaim 
\proclaim{Lemma 1} Let $B$ be a Ferrers board, and let $D$  
be a section of $B$ of width
$d$.  Fix a placement $C$ of $n-d$ rooks in the $n-d$ other columns of the
$n \times n$ grid outside of $D$.  We say a placement of rooks $C^{\prime}$
$extends$ $C$ if all the rooks of $C$ are in $C^{\prime}$. 
Then the minimum value of $\text{mat$(C^{\prime},B)$}$, over
all placements $C^{\prime}$ of $n$ rooks extending $C$, occurs when 
$C^{\prime}$ is $D$-$standard$.  Furthermore, if $E$ is this $D$-$standard$
extension of $C$,
$$
\sum_{C^{\prime}\atop \text{$n$ rooks, extending $C$}}
q^{\text{mat$(C^{\prime},B)$}} = q^{\text{mat$(E,B)$}}[d]!.
$$
\endproclaim \noindent
$Proof:$  By induction on $d$, the case $d=1$ being trivial.
Let $D^{\prime}$ be $D$ minus its left-most column, and call an
extension $C^{\prime}$ {\it $D$ semi-standard} if  
it is $D^{\prime}$-standard.  Also
call the rook in the left-most column of $D$ the {\it left-rook}.
We claim that when we add up $q^{\text{mat}}$ for the $d$ 
semi-standard extensions, we get $q^{\text{mat$(E,B)$}}[d]$.  Lemma 1 will
follow since, if for any of the $d$ choices for the left-rook we let
the $d-1$ remaining rooks in $D$ cycle through their $(d-1)!$
possibilities, by induction they generate an extra $[d-1]!$.

Say the rows left unattacked by the rooks of $C$ are rows
$$
i_s < i_{s-1} < \cdots < i_1, 
$$
(which intersect $B$ within $D$)
and rows 
$$
j_1 < j_{2} < \cdots < j_{d-s}, 
$$
(which do not intersect $B$ within $D$).  Note that $j_1 > i_1$.
We consider what happens
to the statistic $cross$ when we move from one semi-standard placement to
another, in three special cases.

\vskip .1in
\noindent Case $1$:  The left-rook changes from row $j_p$ to row
$j_{p+1}$ for some $1\le p < d-s$.  Consider Figures 5 and 6.
The horizontal line near the middle of the Figure 
is the boundary of $B$, indicated by the letter {\bf B}.  Otherwise,
squares which are counted in the definition of cross are 
indicated by straight lines going through them (ignore 
contributions from rooks outside of $D$ for the moment).  
If, in the definition of $cross$, those squares which
satisfy two of the three conditions were counted twice, then
$cross$ would be the same for Figures 5 and 6, since the
line segments of lengths $P$,$Q$,$R$, and $S$ are mearly
shifted around from one figure to the next.  But there is one
more square in Figure 6 that satisfies two of the three
conditions then there is in Figure 5 (note the circled
intersections; note also that rooks to the left or the right of $D$
will create the same number of intersections in both placements, hence we are
justified in ignoring their contribution when determining how much
$cross$ changes by, and do not need to include them in our figures).  
 Thus $cross$ has decreased by one.

\input epsf
\midinsert
%make all figures 1/2 size
$$\vbox{\epsffile{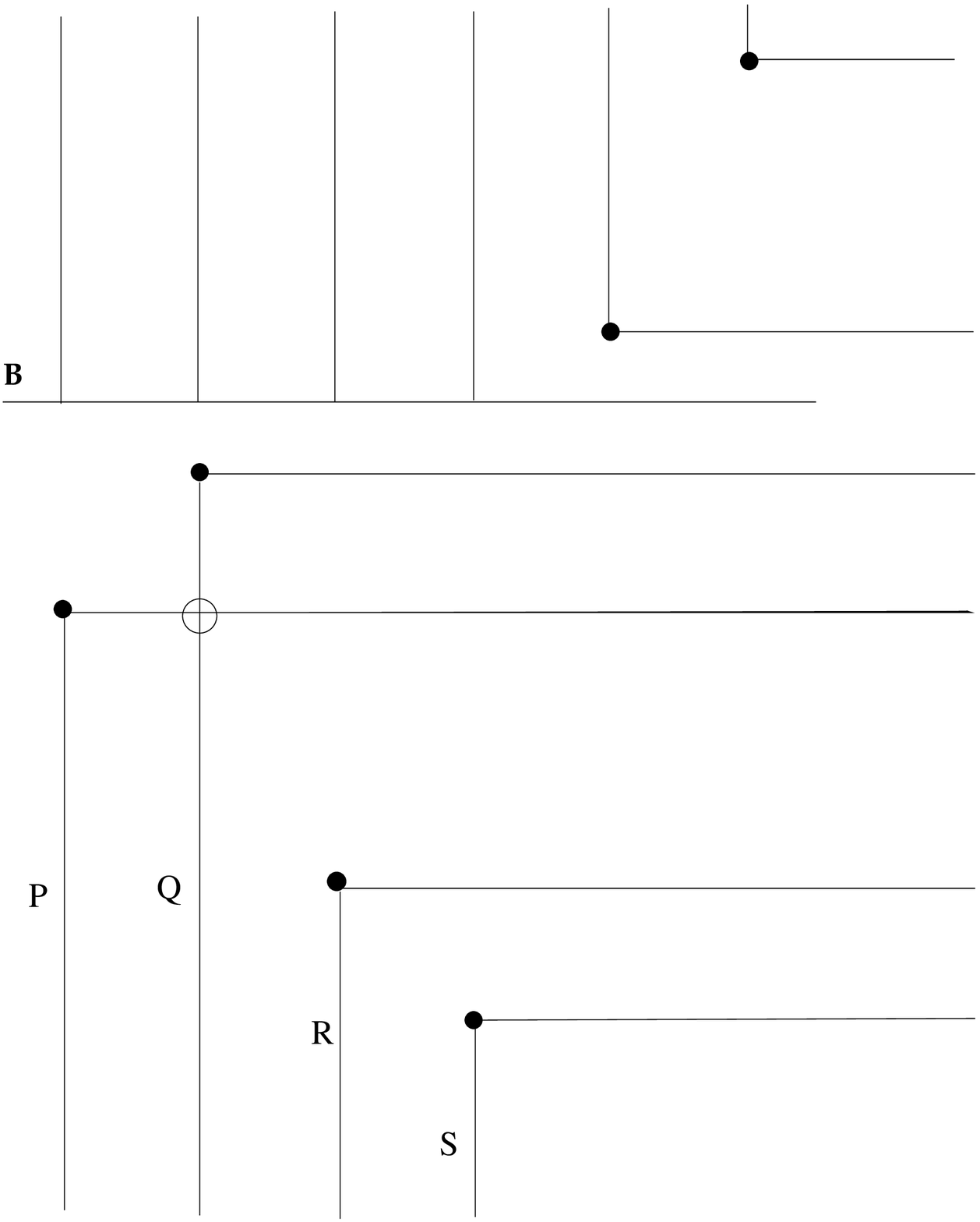}}$$
\botcaption{Figure 5}  A $D$ semi-standard placement of rooks.
\endcaption
\endinsert

\input epsf
\midinsert
%make all figures 1/2 size
$$\vbox{\epsffile{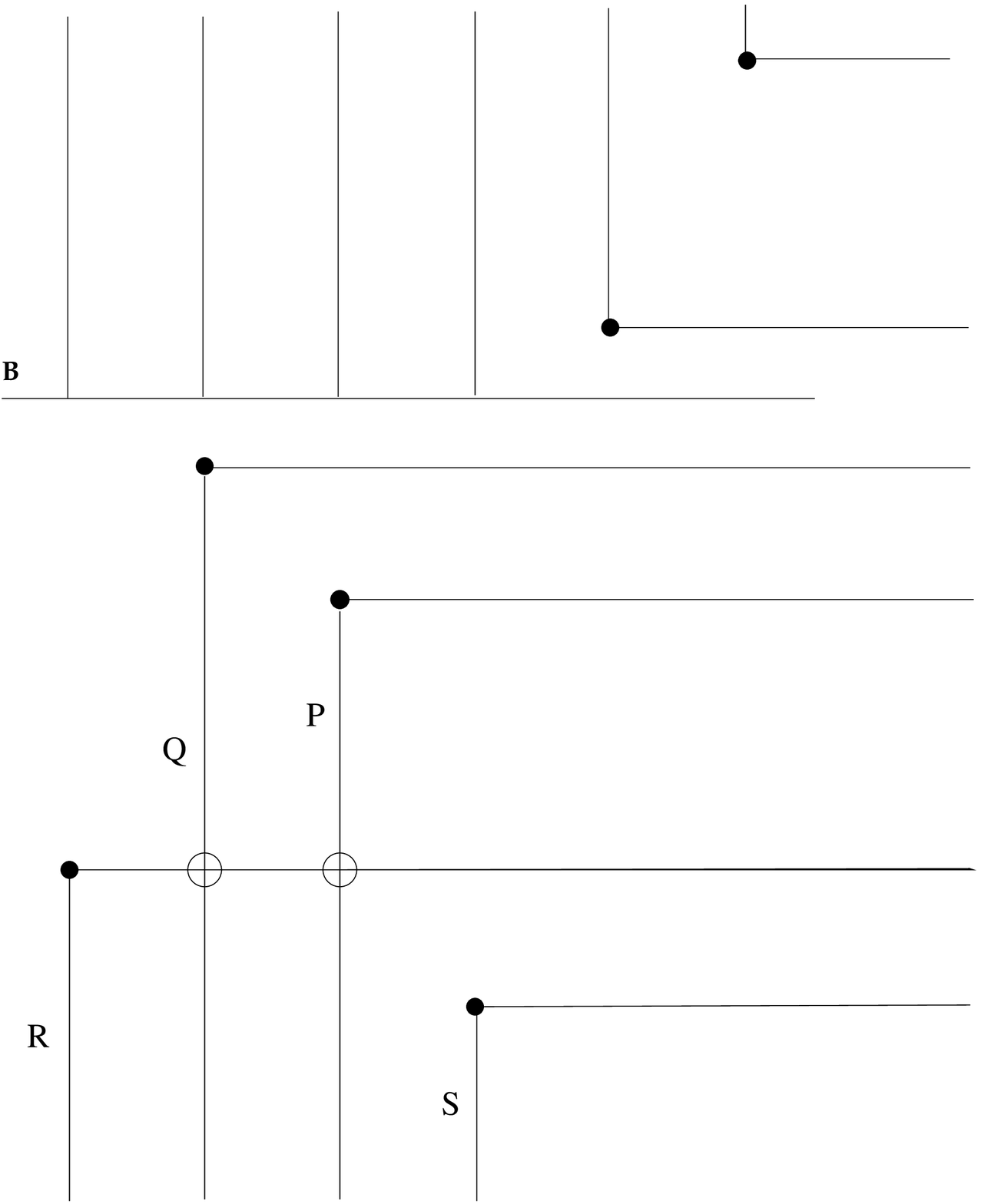}}$$
\botcaption{Figure 6}   
 Another $D$ semi-standard placement of rooks.  The left-rook has moved
 down to row $j_{p+1}$.
\endcaption
\endinsert

\vskip .1in
\noindent Case $2$:  The left-rook changes from row $j_{1}$ to
row $i_1$.  Consider Figures 7 and 8.  As in case 1, we need only
consider the number of squares which satisfy two of the three conditions.
For example, there are $P$ squares below the left-rook in Figure 7,
and also $P$ squares below the rook in column two of
Figure 8.
  Similar remarks apply to $Q$,$R$,$S$,$T$,and $U$.
 Since there are $d-s$ new (circled) intersections,
$cross$ has decreased by $d-s$.

\input epsf
\midinsert
%make all figures 1/2 size
$$\vbox{\epsffile{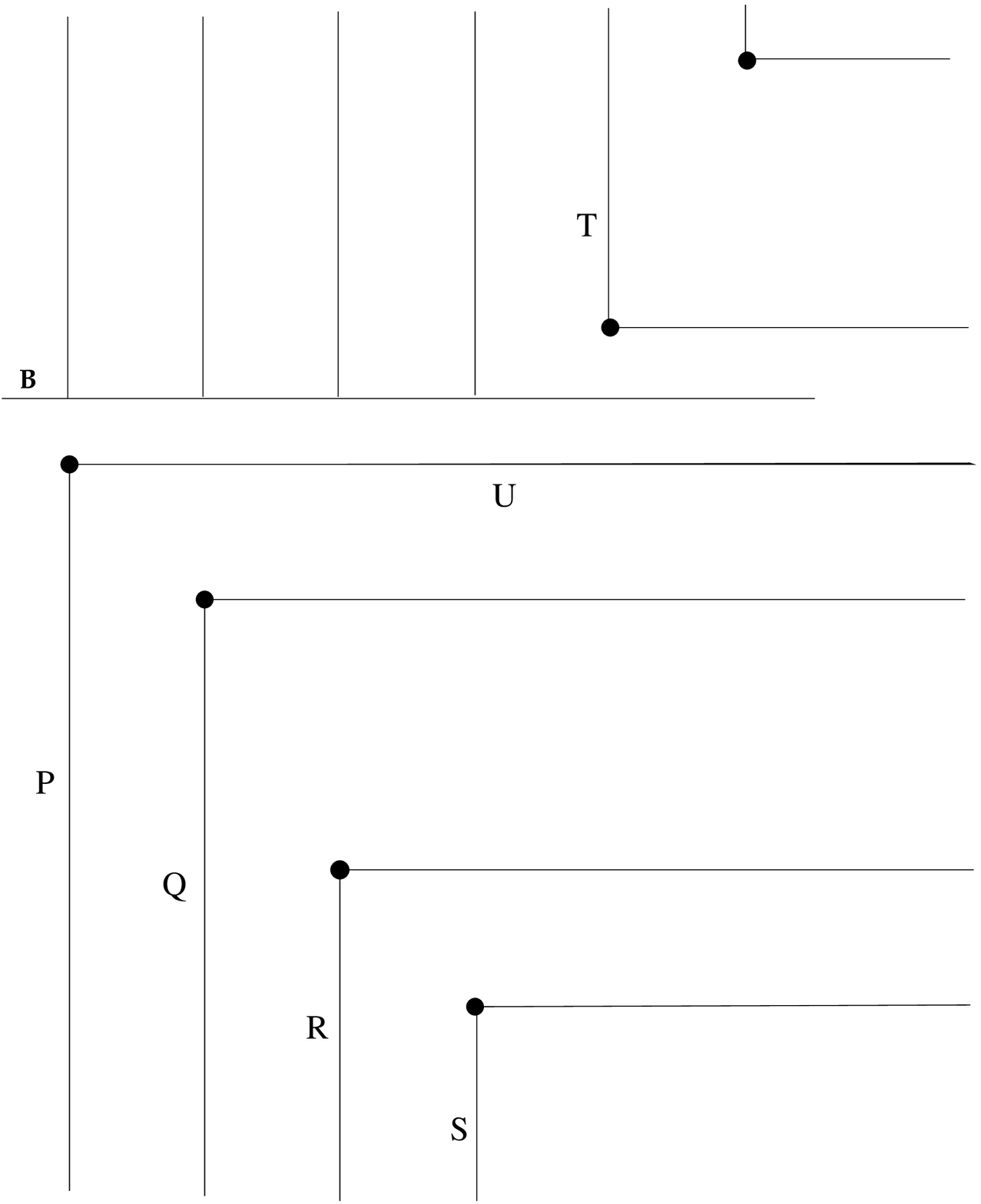}}$$
\botcaption{Figure 7}  A $D$ semi-standard placement of rooks.
\endcaption
\endinsert

\input epsf
\midinsert
%make all figures 1/2 size
$$\vbox{\epsffile{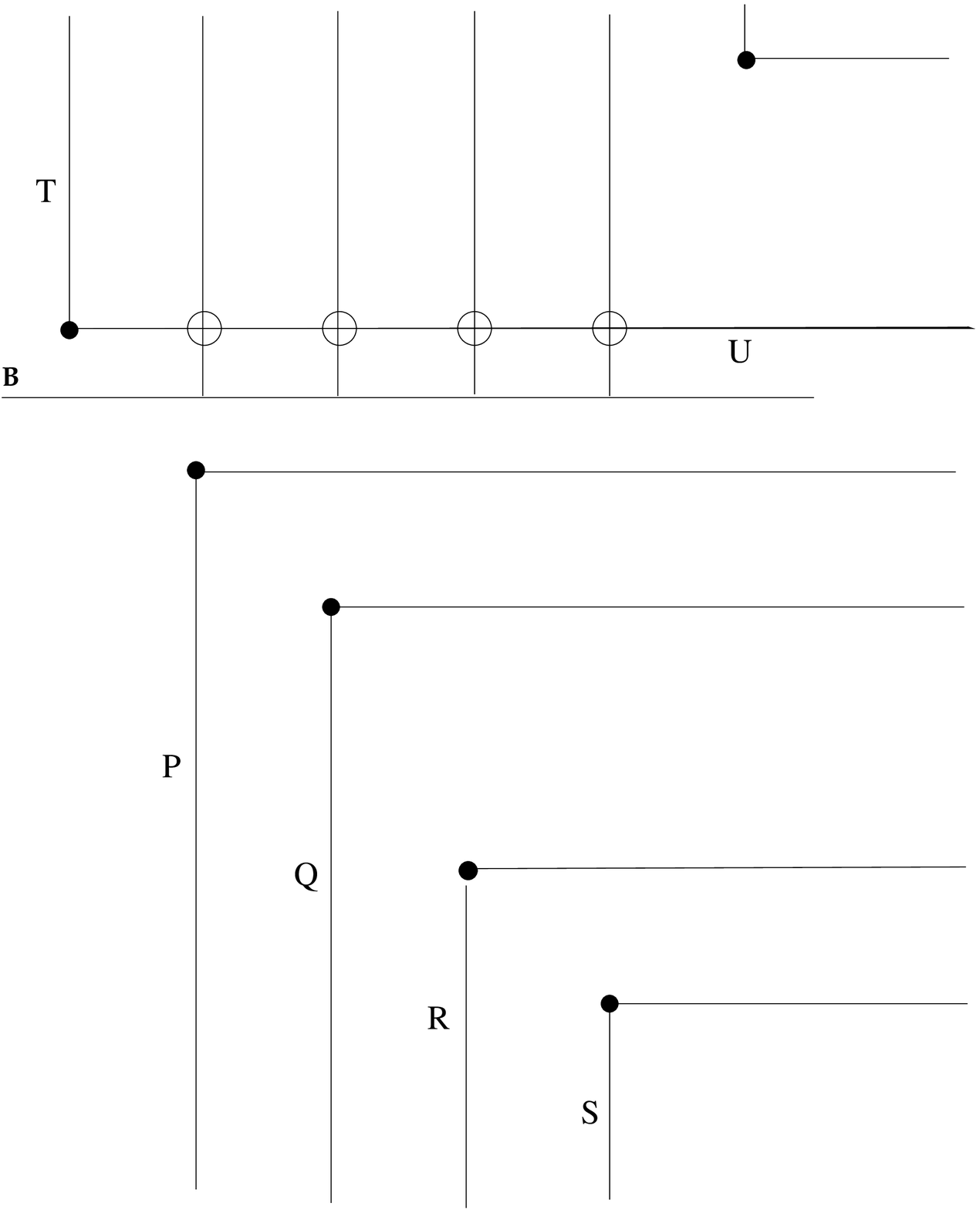}}$$
\botcaption{Figure 8}   
 Another $D$ semi-standard placement of rooks.  The left-rook has moved
 up to row $i_{1}$.
\endcaption
\endinsert

\vskip .1in
\noindent Case $3$:  The left-rook changes from row $i_p$ to row
$i_{p+1}$ for some $1 \le p <s$.  Figures 9 and 10.  As in
case 1, $cross$ decreases by 1.

\input epsf
\midinsert
%make all figures 1/2 size
$$\vbox{\epsffile{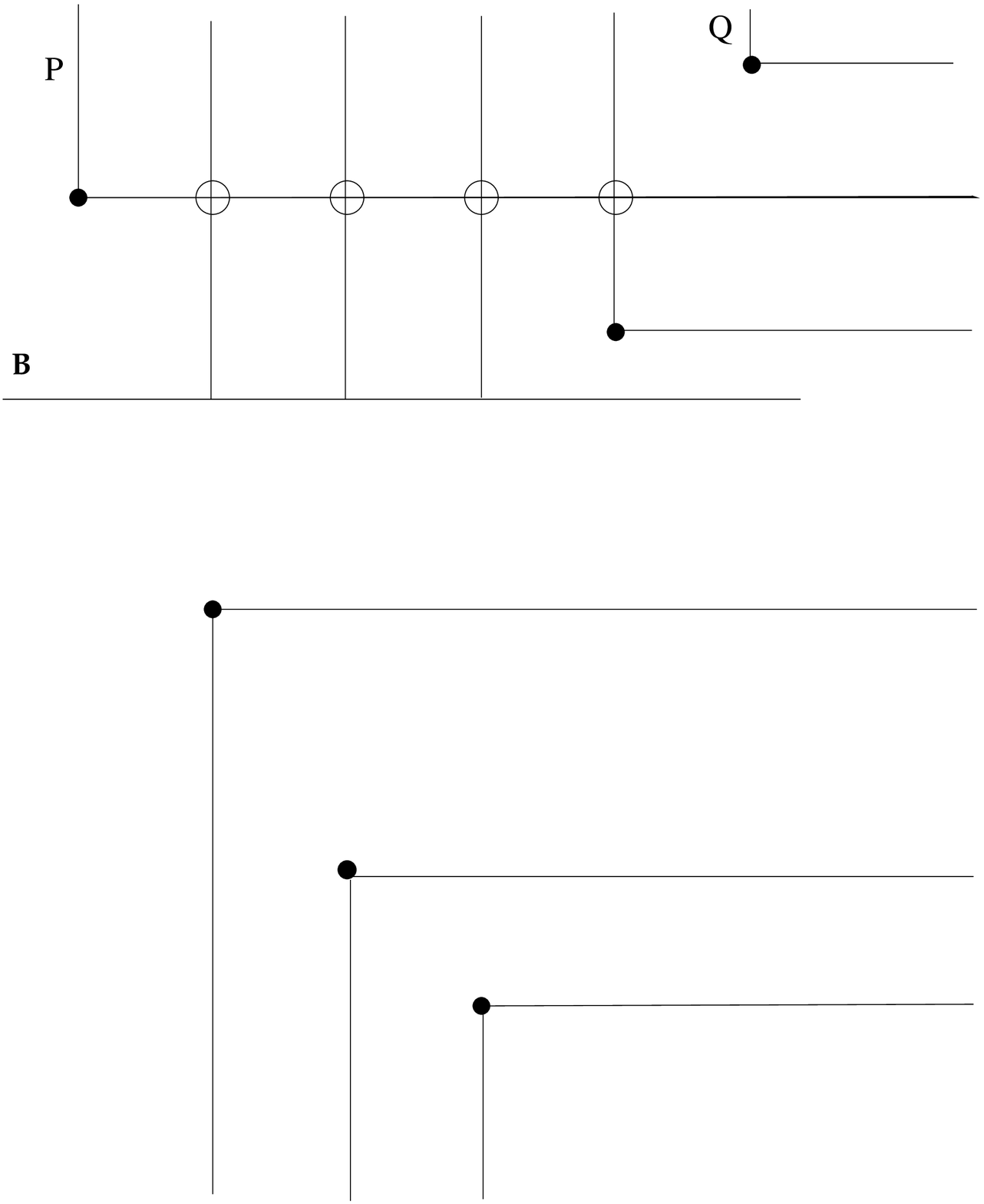}}$$
\botcaption{Figure 9}  A $D$ semi-standard placement of rooks.
\endcaption
\endinsert

\input epsf
\midinsert
%make all figures 1/2 size
$$\vbox{\epsffile{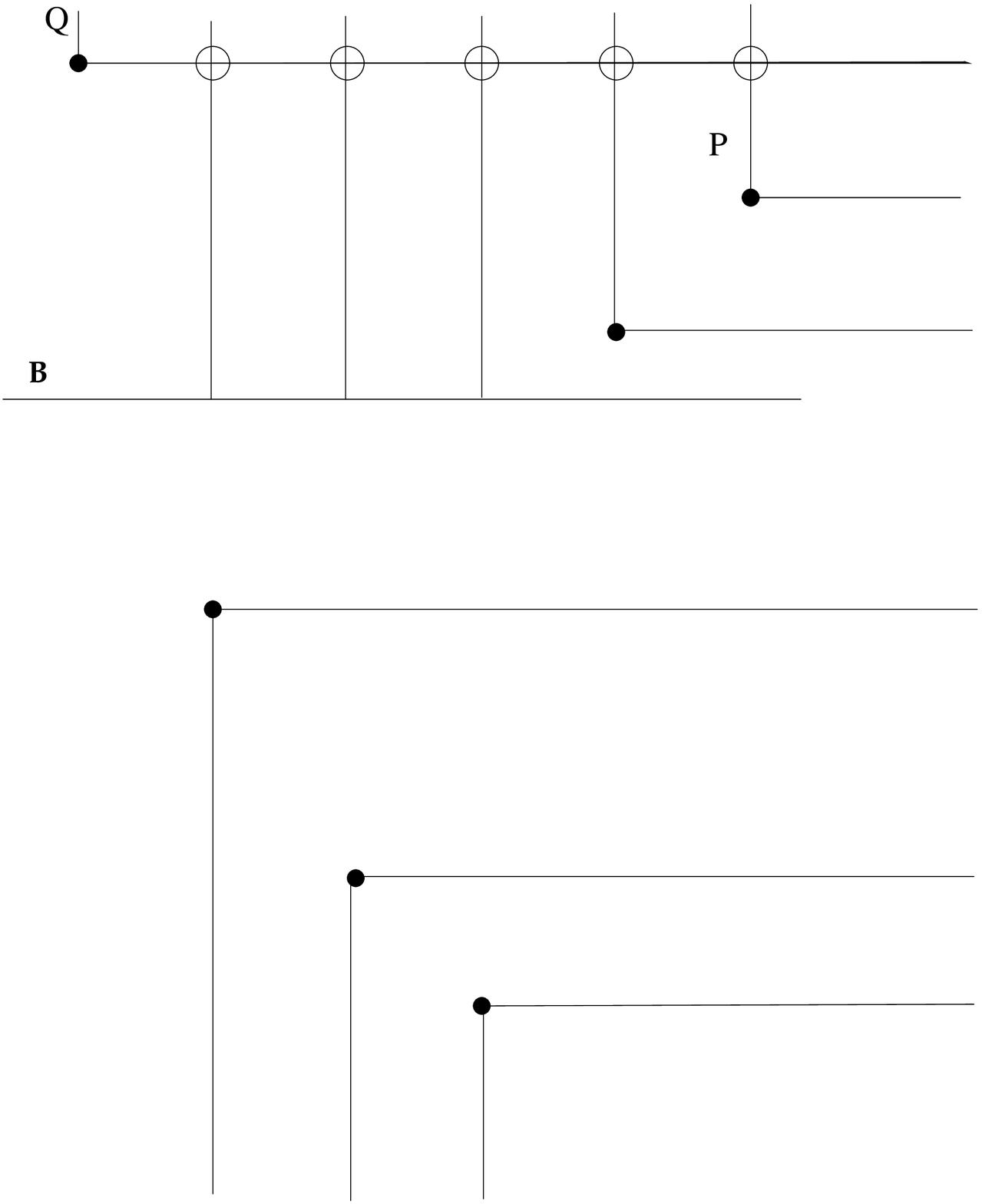}}$$
\botcaption{Figure 10}   
 Another $D$ semi-standard placement of rooks.  The left-rook has moved
 up to row $i_{p+1}$.
\endcaption
\endinsert

\vskip .1in

Combining cases 1,2, and 3, we see that as the left-rook cycles through
rows 
$$
j_1,j_2,\ldots, j_{d-s},i_1,i_2,\ldots ,i_s,
$$
$cross$ decreases
by one each time, hence $mat$ increases by one each time.  Thus
$$
\sum_{\text{semi-standard $C^{\prime}$} \atop 
\text{extending $C$}}
q^{\text{mat$(C^{\prime},B)$}}= q^{\text{mat$(E,B)$}}[d]
$$
and Lemma 1 follows by induction. $\qquad \blacksquare$

Let $B$ be the Ferrers board of Fig. 11, with
$\bold d=(d_1,\ldots ,d_t) \in \Bbb P^t$ a vector satisfying $\sum_i d_i=n$.
For technical reasons we allow the $h_i \in \Bbb N$ (thus there are in
general several different choices 
for $t$, $\bold h$, and $\bold d$ which represent
the same board).
  A placement $C$ of $n$ rooks on squares 
$(1,\tau_1),\ldots ,(n,\tau_n)$ can be converted into a multiset
permutation $\sigma \in M(\bold d)$ by first 
forming a sequence $S$ whose $i$th element
is $\tau_i$, then replacing numbers 1 through $d_1$ of $S$, wherever they
occur, by all
$1$'s, numbers $d_1+1$ through $d_2$ of $S$ by all $2$'s, etc.  We call the
$\prod_id_i!$ placements that get mapped to $\sigma \in M(\bold d)$ 
the $graphs$ of
$\sigma$. 
By restricting our attention to
$B$-standard placements, we have a bijection between elements
$\sigma \in M(\bold d)$ and $B$-standard placements $C(\sigma)$.  We call $C(\sigma)$ 
the $B$-$standard$ $graph$ of $\sigma$.

\proclaim {Definition 4} Let 
$B=B(h_1,d_1;h_2,d_2;\ldots ;h_t,d_t)$ be the Ferrers board of
Fig. 11, where $\sum_i d_i=n$.
For a given permutation $\sigma \in M((d_1,d_2,\ldots ,d_t))$, let
$$
\text{mat}(\sigma,B):=\text{mat}(C(\sigma),B)
$$ 
where $C(\sigma)$ is the $B$-standard graph of  
$\sigma$.
\endproclaim

\input epsf
\midinsert
%make all figures 1/2 size
$$\vbox{\epsffile{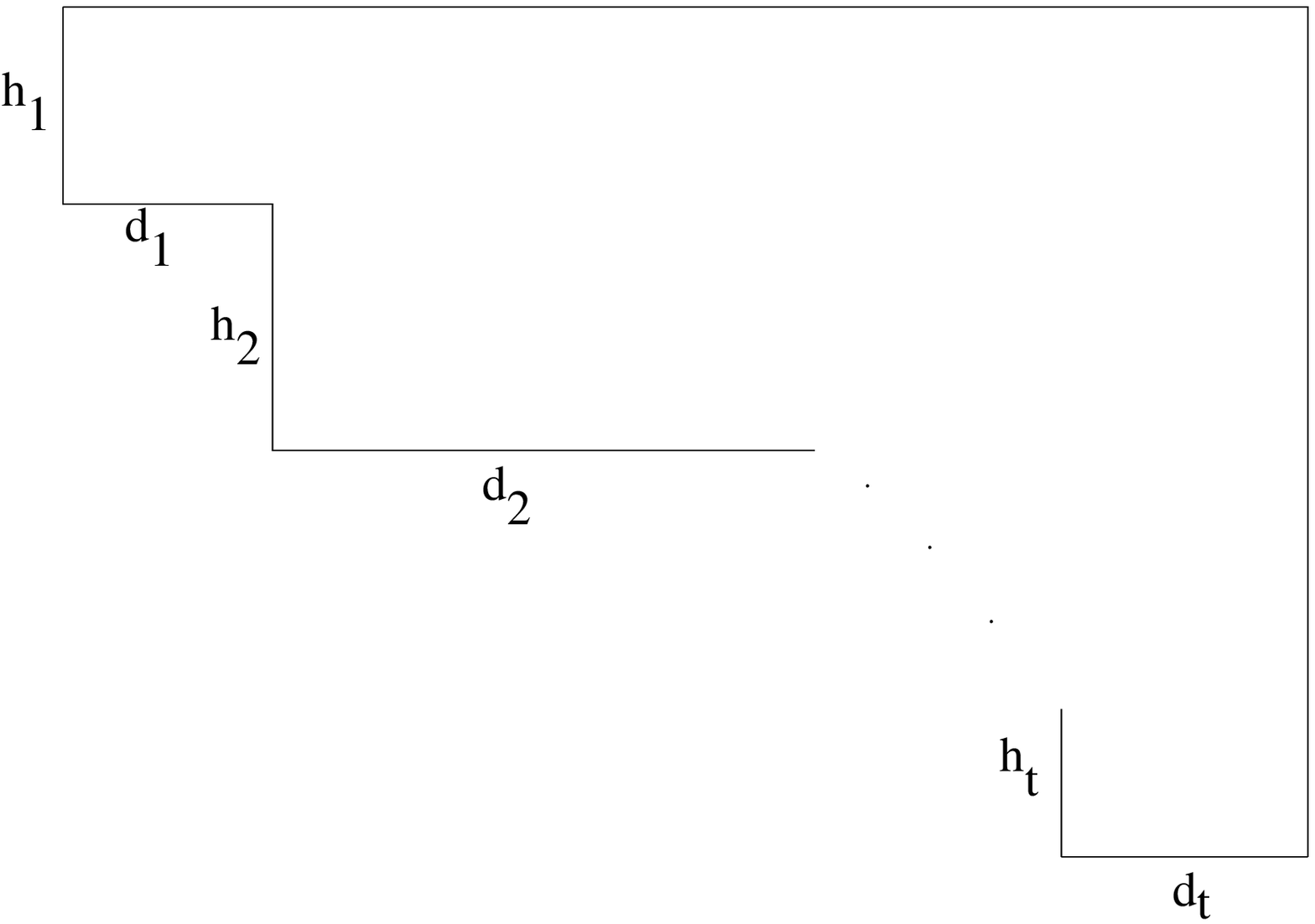}}$$
\botcaption{Figure 11}  The Ferrers board $B(h_1,d_1;\ldots ;h_t,d_t)$,
where $d_i \in \Bbb P, h_i \in \Bbb N \text{ for } 1\le i \le t$.
The first $d_1$ columns have height $h_1$, the next $d_2$ have height
$h_1+h_2$, etc. 
\endcaption
\endinsert

By iterating Lemma 1 and using (11) we now have 
\proclaim {Theorem 3}  Let $B:=B(h_1,d_1;\ldots ;h_t,d_t)$ 
be the Ferrers board of Fig. 11, with $\sum_i d_i=n$.
 Then $\text{mat}(B)$ is
multiset Mahonian, i.e. 
$$
\sum_{\sigma \in M(\bold d)} q^{\text{mat}(\sigma,B)} =
\left [ \matrix & n & \\ d_1, & d_2, & \ldots &,d_t \endmatrix
\right ].
$$
\endproclaim

\proclaim {Definition 5} Let $B$ be a Ferrers board, and let $D$ 
be a section of $B$ of width $d$.  Let $C$ be a placement of 
$n$ rooks on the
$n \times n$ grid, with $s$ of the $d$ rooks in $D$ on $B$.  
We say $C$ is $D$-$regular$ 
if both of the following hold.
\vskip .1in
\noindent 
1) The $s$ rooks on $B$ and in $D$ are in the $s$ left-most columns of $D$.
   Furthermore, these $s$ rooks are in ``ascending" order; if two of these
   rooks occupy squares $(i,j)$ and $(k,l)$, with $i< k$, then $j> l$.
\vskip .1in

\noindent
2) The $d-s$ remaining rooks in $D$ and off $B$ are also in ascending order.
\vskip .1in
\noindent
Call a placement of $n$ rooks on the
$n \times n$ grid $B$-$regular$ if it is $D$-$regular$ for all possible sections
$D$ of $B$.  See Fig. 12.
\endproclaim

\input epsf
\midinsert
%make all figures 1/2 size
$$\vbox{\epsffile{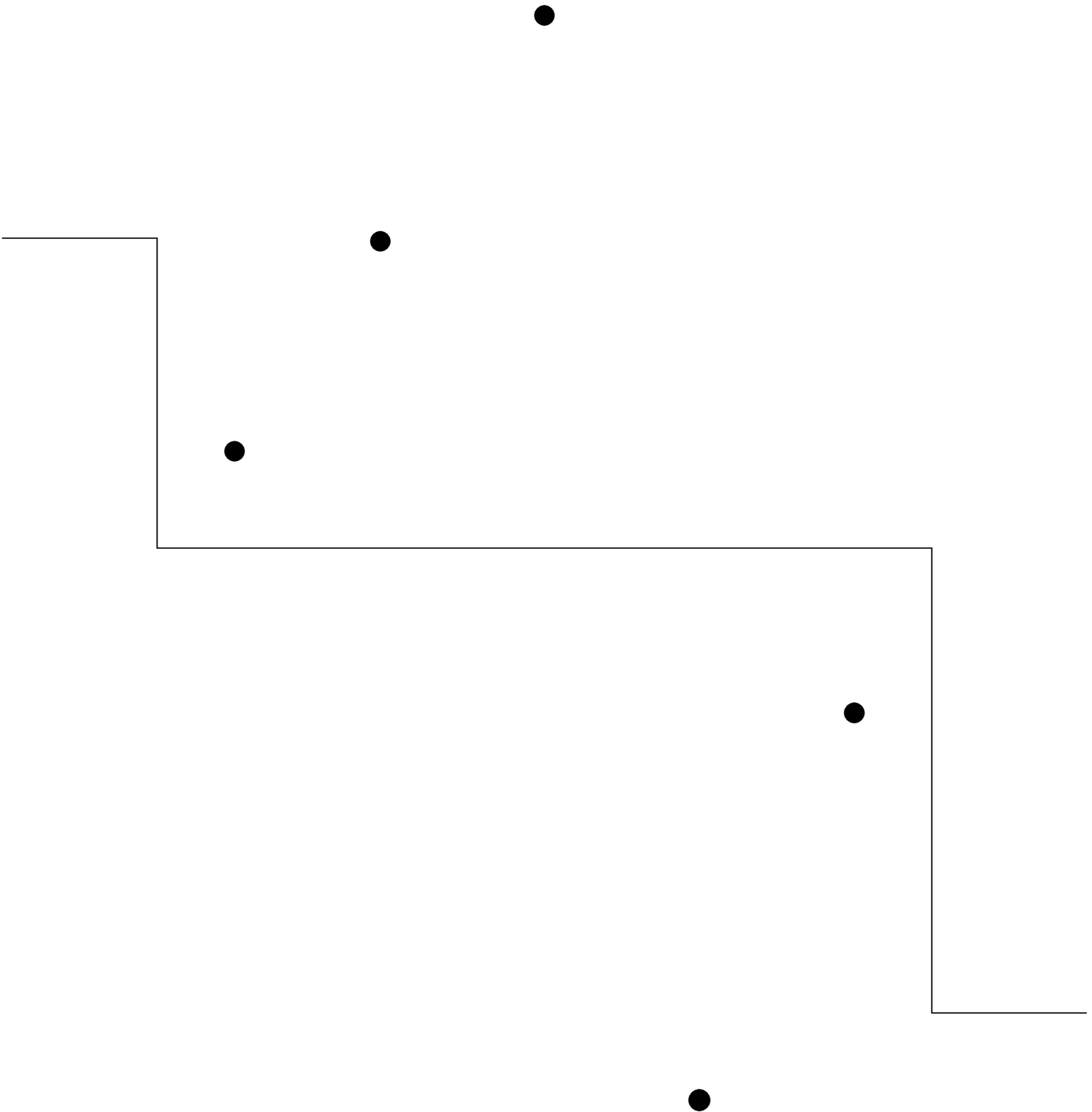}}$$
\botcaption{Figure 12}  A section $D$ of width $5$ and a
$D$-regular placement of rooks (rooks outside of $D$ are not
pictured).
\endcaption
\endinsert

By arguments identical to those given for the statistic $mat$, one can prove
the following results for $\xi$. 
\proclaim{Lemma 2} Let $B$ be a Ferrers board, and let $D$ be a section of width
$d$.  Fix a placement $C$ of $n-d$ rooks in the $n-d$ other columns of the
$n \times n$ grid outside of $D$.  
Then the minimum value of $\xi (C^{\prime},B)$, over
all placements $C^{\prime}$ of $n$ rooks extending $C$, occurs when 
$C^{\prime}$ is $D$-$regular$.  Furthermore, if $E$ is this $D$-$regular$
extension of $C$,
$$
\sum_{C^{\prime}\atop \text{$n$ rooks, extending $C$}}
q^{\xi (C^{\prime},B)} = q^{\xi (E,B)}[d]!.
$$
\endproclaim 

\proclaim {Definition 6} Let 
$B=B(h_1,d_1;h_2,d_2;\ldots ;h_t,d_t)$ be the Ferrers board of
Fig. 11, where $\sum_i d_i=n$.
For a given permutation $\sigma \in M((d_1,d_2,\ldots ,d_t))$, let
$$
\xi(\sigma,B):=\xi (P(\sigma),B),
$$ 
where $P(\sigma)$ is the graph of $\sigma$ which is $B$-regular.
\endproclaim

\proclaim {Theorem 4}  Let $B:=B(h_1,d_1;\ldots ;h_t,d_t)$ 
be the Ferrers board of Fig. 11, with $\sum_i d_i=n$.
 Then $\xi(B)$ is
multiset Mahonian, i.e. 
$$
\sum_{\sigma \in M(\bold d)} q^{\xi(\sigma,B)} =
\left [ \matrix & n & \\ d_1, & d_2, & \ldots & ,d_t \endmatrix
\right ].
$$
\endproclaim

\head 4. Euler-Mahonian Statistics \endhead

\proclaim {Definition 7}  Let 
$B(n)$ denote the triangular board of side $n-1$ consisting of all
squares $(i,j)$ with $1\le i< j\le n$.
\endproclaim

In [Ha1] it was shown that
$$
\sum_{\sigma \in S_n \atop \text{des}(\sigma)=k}q^{\text{maj}(\sigma)}
=q^{nk-{n \choose 2}}T_{k}(B(n)) \eqno (12)
$$
and also that
$$
\sum_{\sigma \in S_n \atop \text{des}(\sigma)=k}q^{\text{maj}(\sigma)}
=T_{n-k-1}(B(n)). \eqno (13)
$$
For any Ferrers board $B(c_1,\ldots,c_n)$, let $B^c$ be the complementary
board with column heights $n-c_n,n-c_{n-1},\ldots ,n-c_1$.  Dworkin proved
a ``reciprocity" theorem for $T_k$, namely [Dwo]
$$
T_k(B,q^{-1})=q^{-{n \choose 2}}T_{n-k}(B^c,q). \eqno (14)
$$
  Combining (13), (14), and a symmetry 
property of the $T_k$ (Theorem 6) which we prove in section 5, we get
$$
\sum_{\sigma \in S_n \atop \text{des}(\sigma)=k}q^{\text{maj}(\sigma)}
=q^{nk-{n \choose 2}}T_{k+1}(B(n)^c,q). \eqno(15)
$$
Similarly,
(12) and (14) imply 
$$
\sum_{\sigma \in S_n \atop \text{des}(\sigma)=k}q^{\text{maj}(\sigma)}
=T_{n-k}(B(n)^c,q), \eqno (16)
$$
a fact we will use later.

There is a straightforward way, used by Riordan and Kaplansky [KaRi], to
identify a permutation $\pi=\pi_1 \pi_2 \cdots \pi_n \in S_n$ having
$k$ descents, with a placement $F(\pi)$ of $n$ rooks 
on the $n \times n$ grid with 
$k$ rooks on $B(n)$.  If $\pi_{j_1}=1$, let $y_1$ be the cycle
$(\pi_1 \pi_2 \cdots \pi_{j_1})$.   If $\alpha$ is the smallest integer 
not contained in $y_1$, and $\pi_{j_2}=\alpha$, let $y_2$ be the cycle
$(\pi_{j_1+1}\pi_{j_1+2}\cdots \pi_{j_2})$, etc.
 Now let $F(\pi)$ be the placement having a rook on $(i,j)$ if and
 only if $i$ and $j$ are in the same cycle $y_p$ for some $p$, with
 $i$ immediately following $j$.
  Call $F(\pi)$ the $descent$ $graph$ of 
$\pi$.  
 For example, if $\pi=3521647$, the $y_i$ are the cycles
$(3521)$, $(64)$, and $(7)$, and its descent graph is illustrated in 
in Fig. 13.

\input epsf
\midinsert
%make all figures 1/4 size
$$\vbox{\epsffile{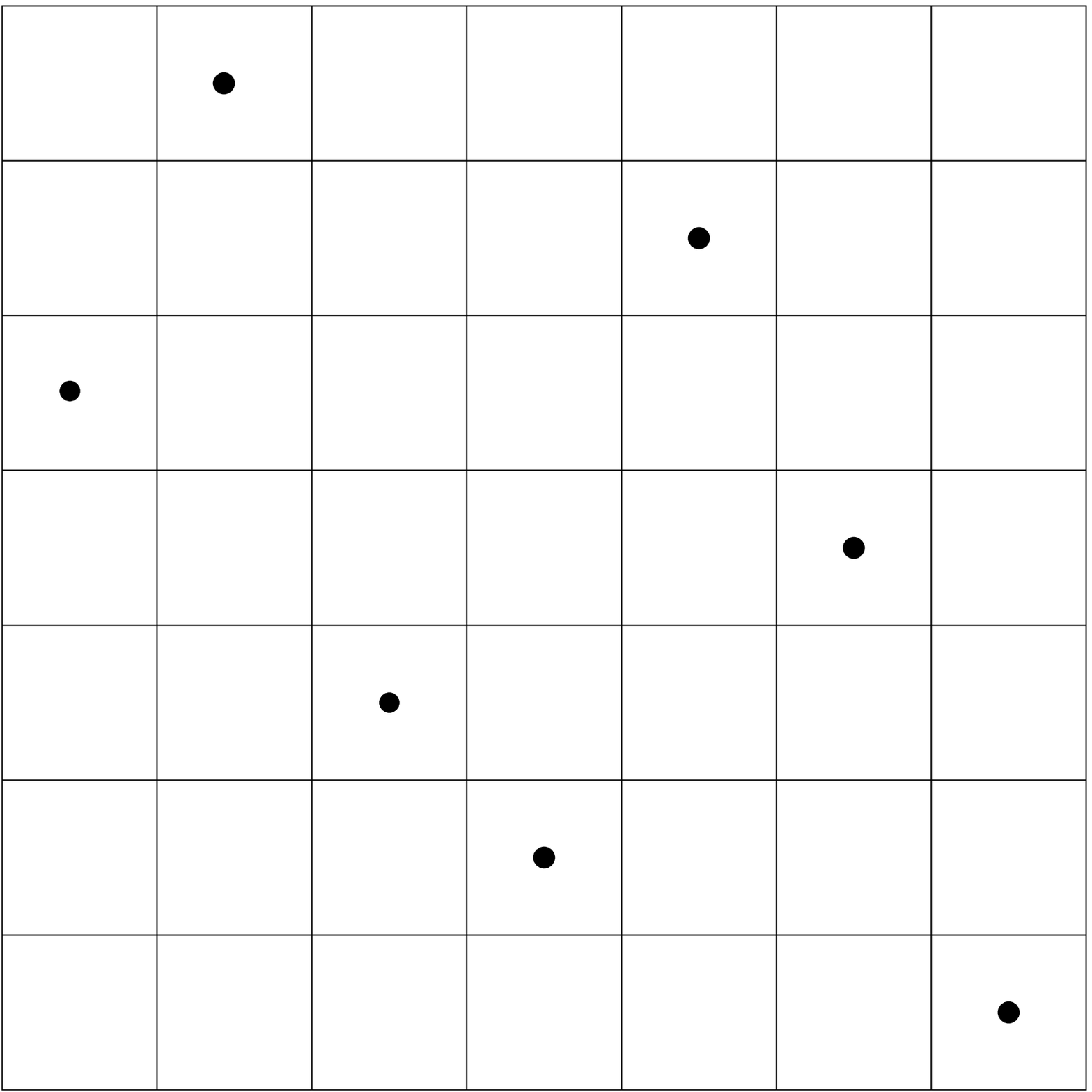}}$$
\botcaption{Figure 13} The descent graph of 
the permutation $3521647$, whose associated cycles are
$(3521)$, $(64)$, and $(7)$.
There are three rooks on $B(n)$, corresponding to the three descents of
$3521647$.
\endcaption
\endinsert

It follows from (12) that if we
define 
$$
\text{stat1}(\pi):= n\text{des}(\pi) - {n \choose 2}+\text{mat}(F(\pi),B(n)) 
\eqno (17)
$$
$$
= n^2-cross(F(\pi),B(n)), 
$$
then $(des,stat1)$ is jointly distributed with $(des,maj)$. 
To get another such pair $(des,stat2)$ we
can
reflect the board about the cross diagonal, i.e. relabel square
$(i,j)$ as square $(n-j+1,n-i+1)$, which gives us a new
rook placement $F^{\prime}(\pi)$ with the same number of rooks on $B(n)$.
For example, if we reflect the placement in Fig. 13 we get the
descent graph of $1425763$.
This placement will have a different value of $mat$, which we
can then use to define $stat2$ as in (17) above.

If we reverse a permutation $\pi$ with 
$k$ descents, we 
get a new permutation with $n-k+1$ descents, say $\beta(\pi):=\pi_n \pi_{n-1}\cdots \pi_1$. 
  By (13), 
if we let $\text{stat3}(\pi)=\text{mat}(F(\beta(\pi)),B(n))$ 
, we have an Euler-Mahonian pair $(des,stat3)$.  We can also get
another pair $(des,stat4)$ by reflection.  So far we have four statistics for both $\xi$
and $mat$ which, when combined with $des$, form an Euler-Mahonian
pair.  For each of these statistics $stat$ we
can get another Euler-Mahonian pair by forming $(des,nk-stat)$ (if we let
$\zeta_i (\pi):=n-\pi _{n-i+1}+1$, then $\text{des}(\zeta (\pi))=\text{des}(\pi)$
and $\text{maj}(\zeta (\pi))=n\text{des}(\pi)-\text{maj}(\pi)$, hence the LHS
of (13) is symmetric about $q^{nk/2}$).  Thus both $\xi$ and
$mat$ each induce a family of $8$ pairs.  Table $1$ lists the two families for
three sample permutations.

\input epsf
\midinsert
%make all figures normal size
$$\hss\vbox{\epsffile{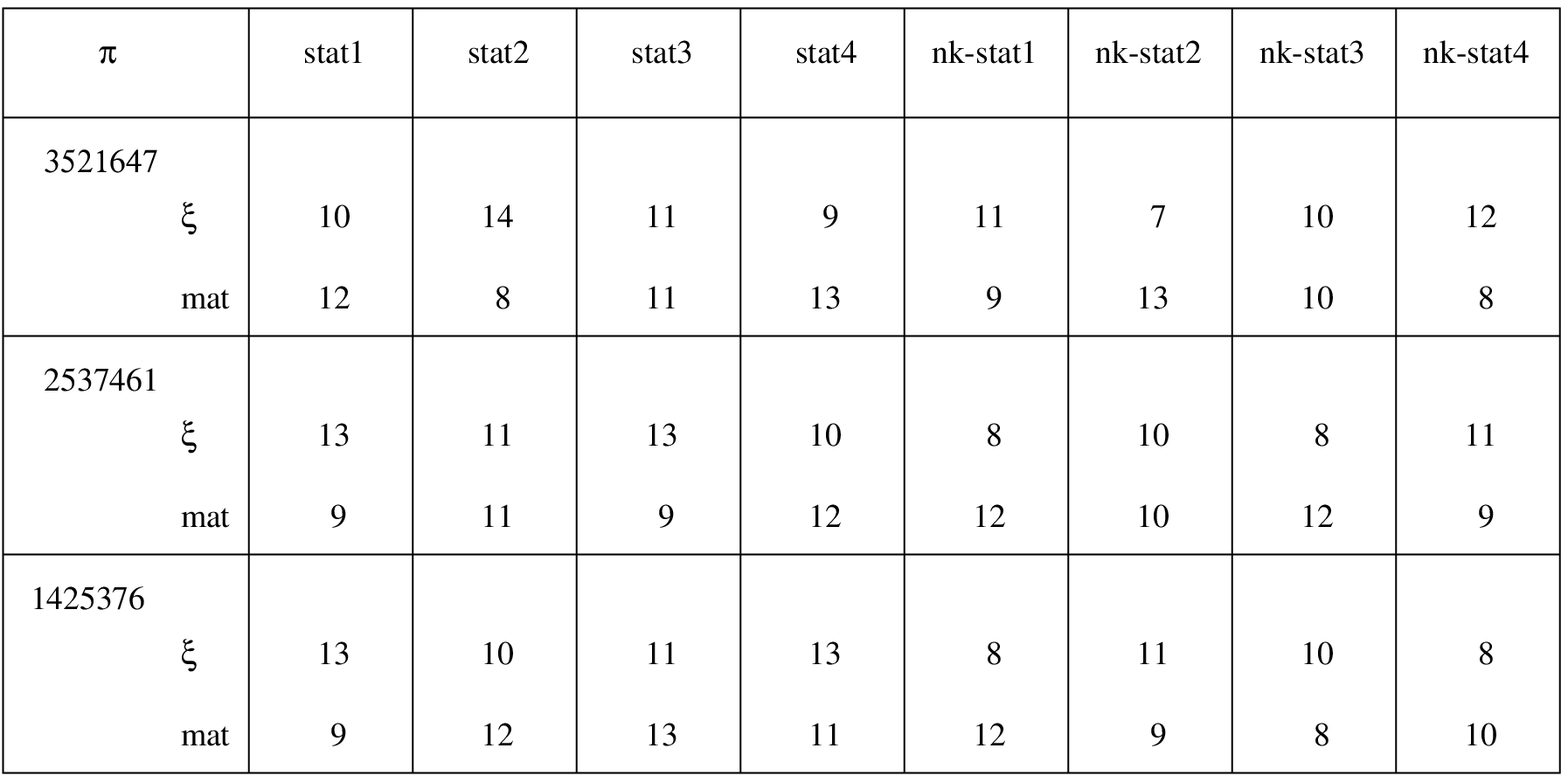}}\hss$$
\topcaption{Table 1} 
\endcaption
\endinsert
 
An examination of Table $1$ shows that none of the 
sixteen pairs equal each other for all $\pi$.  Hence 
the $\xi$ and $mat$ families are
fundamentally different,
at least with respect to
the simple transformations we have considered here.  In addition none of the
sixteen pairs are equal to $maj$, and are also unequal to the statistic $mak$ as
described in [CSZ] (the pair $(des,mak)$ is known to be Euler-Mahonian [FoZe]).

We can also make use of (15) and (16) to try and generate other Euler-Mahonian
pairs.  However, examples indicate that the pairs arrived at in this manner are
rearrangements of the
sixteen pairs above.
%, and we content ourselves here with simply listing one of these
%presumed relationships.  
%
%If we For each $1\le j\le 4$, abbreviate the statistic
%corresponding to the $\xi$ version of $statj$ in column $j$ of Table 1 by
%$\xi_j$, and the $mat$ version by $m_j$.  Furthermore, abbreviate the

\input epsf
\midinsert
%make all figures 1/2 size
$$\vbox{\epsffile{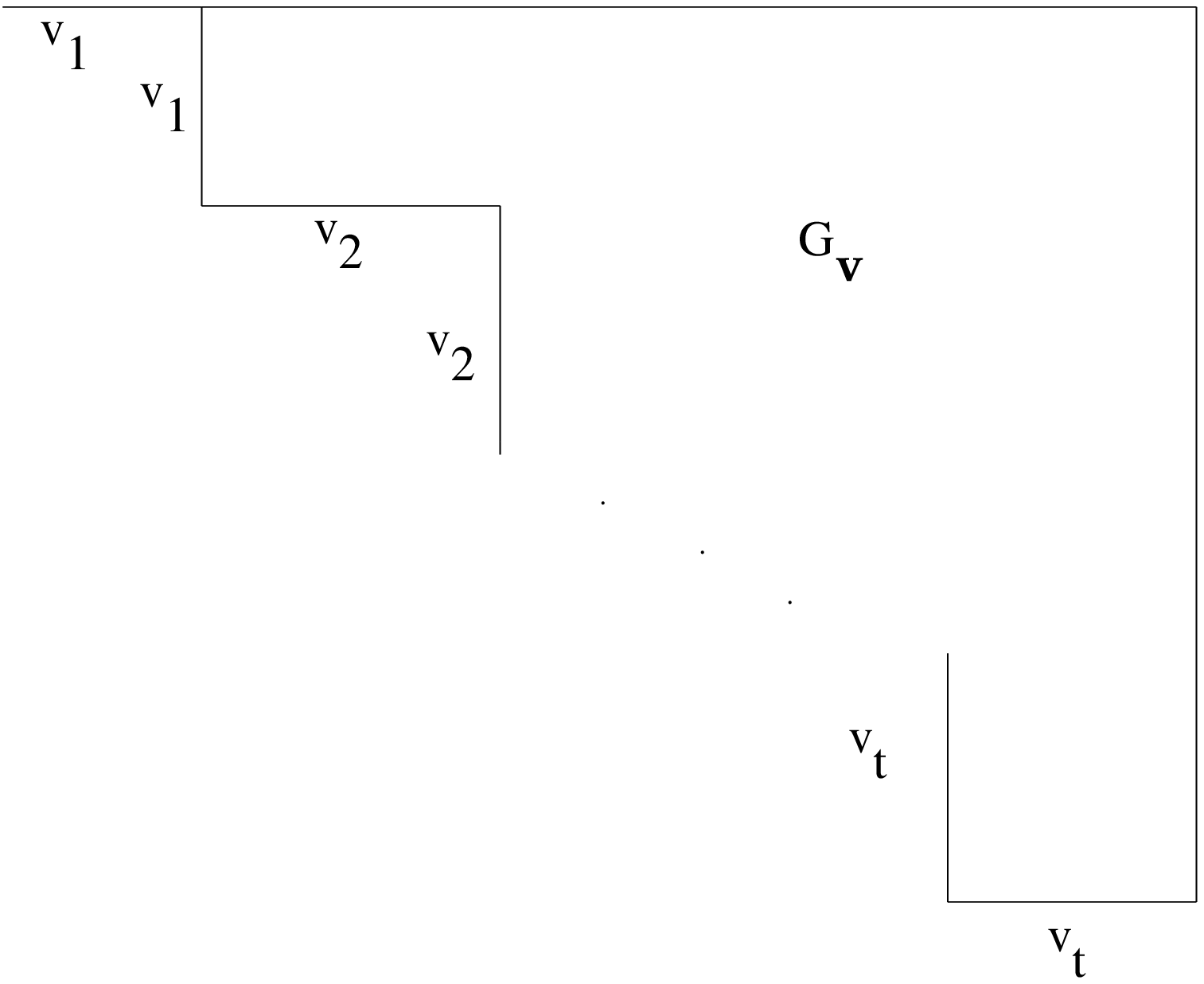}}$$
\botcaption{Figure 14}  The Ferrers board $G_{\bold v}$.  The first
$v_1$ columns are empty, the next $v_2$ have height $v_1$, etc.
\endcaption
\endinsert

Some of our Euler-Mahonian statistics can easily be rephrased as a
multiset solution to (4).  We utilize the following generalization of (12)
[Hag1,p.118];
$$
\sum_{\pi \in M(\bold v) \atop \text{des}(\pi)=k} 
q^{\text{maj}(\pi)} =
\frac{T_{k}(G_{\bold v})}{\prod_{i=1}^t [v_i]!}
q^{nk-\text{Area}(G_{\bold v})}, \eqno(18)
$$
where $G_{\bold v}$ is the board of Fig. 14.
For $\sigma \in M(\bold v)$, it is easy to see that
the number of rooks on $G_{\bold v}$ in any of the $\prod_iv_i!$
graphs of $\sigma$ is $\text{exc}(\sigma)$.  Using this Lemma 1 and (18) imply
$$
\sum_{\text{$G_{\bold v}$-standard placements $C$} \atop 
\text{$k$ rooks on $G_{\bold v}$}}
q^{nk-\text{Area}(G_{\bold v})
+\text{mat}(C,G_{\bold v})} =
\sum_{\sigma \in M(\bold v)\atop \text{des}(\sigma)=k}
q^{\text{maj}(\sigma)} 
$$
$$
=
\sum_{\sigma \in M(\bold v)\atop \text{exc}(\sigma)=k}
q^{\text{stat}5(\sigma)}, \eqno(19)
$$
where $\text{stat}5(\sigma):=nk - \text{Area}(G_{\bold v})
+\text{mat}(C(\sigma),G_{\bold v})$,
with $C(\sigma)$ the $G_{\bold v}$-standard graph of  
$\sigma$.
We also have
$$
\sum_{\sigma \in M(\bold v)\atop \text{exc}(\sigma)=k}
q^{\text{stat}6(\sigma)} 
=
\sum_{\sigma \in M(\bold v)\atop \text{des}(\sigma)=k}
q^{\text{maj}(\sigma)} 
$$
where $\text{stat}6(\sigma):=\xi(P(\sigma),G_{\bold v})
+nk -\text{Area}(G_{\bold v})$, 
with $P(\sigma)$ the $G_{\bold v}$-regular graph 
of $\sigma$.

If we reflect the board $G_{\bold v}$, the $G_{\bold v}$-standard graph of 
$\sigma \in M({\bold v})$ gets sent to one of the $\prod_i v_i!$ graphs of
some $\overline{\sigma} \in M((v_t,v_{t-1},\ldots ,v_1))$, with
$\text{exc}(\sigma)=\text{exc}(\overline{\sigma})$.  Thus we have the
identity
$$
\sum_{\sigma \in M((v_t,\ldots ,v_1)) \atop \text{exc}(\sigma)=k} 
q^{\text{stat}7(\sigma)} =
\sum_{\pi \in M((v_1,\ldots ,v_t)) \atop \text{des}(\pi)=k} 
q^{\text{maj}(\pi)},  \eqno (20)  
$$
where $\text{stat}7(\sigma)=nk-\text{Area}(G_{\bold v})+
\text{mat}(C({\overline \sigma}),G_{\bold v})$, with ${\overline \sigma}$
and $\sigma$ the 
reflected images of each other.  MacMahon showed the RHS of (20) is invariant
under any permutation of the coordinates of $\bold v$ (the RHS of
(2) reduces to $\prod_i \left [ \matrix x+v_i-1 \\ v_i \endmatrix
\right ][v_i]!$ when $B=G_{\bold v}$, which implies the LHS of
(2) and hence the $R_k$ are so invariant, and thus 
also the $T_k$ by (3)).  This gives
$$
\sum_{\sigma \in M((v_t,\ldots ,v_1)) \atop \text{exc}(\sigma)=k} 
q^{\text{stat}7(\sigma)} =
\sum_{\pi \in M((v_t,\ldots ,v_1)) \atop \text{des}(\pi)=k} 
q^{\text{maj}(\pi)},    
$$
a proper multiset Euler-Mahonian pair.
Clearly we can also replace $\text{stat}7$ above by
$\xi (P({\overline \sigma}),G_{\bold v})$
$+nk-\text{Area}(G_{\bold v})$, 
where $P({\overline \sigma})$ is the $G_{\bold v}$-regular
graph of ${\overline \sigma}$.

One could also generate other identities by applying reciprocity to (18),
but
instead of working with excedences, we would need to work with $rises$,
a rise being a value of $i$ such that $\sigma_i \ge f_i(\bold v)$.

It doesn't seem to be as easy to obtain new statistics by
reversing the string when working with multiset
permutations since if such a permutation has $k$ excedences (or $k$ descents), the number
of excedences (or descents) of the reversed string is unpredictable.

In the late 1980's M. Denert introduced an interesting permutation statistic 
which arose during her research into algebraic number theory.  She
conjectured that this statistic was Euler-Mahonian when paired with $exc$.
 Her conjecture was
proven by Foata and Zeilberger [FoZe], who named her statistic ``Denert's statistic",
denoted by $den$.  
We now show that $den$ is part of the
$\xi$ family.  

It will prove convenient to work with $B(n)^c$ and (16).  Let
$\sigma \in S_n$, and let $C(\sigma)^T$ denote the transpose of the graph of $\sigma$
(which is the graph of $\sigma^{-1}$).  For each rook in $C(\sigma)^T$,
put an $x$ on the grid in all the squares to the right and in the row.
For each rook off $B(n)^c$, put a circle in the squares below and in the
column, and also 
in the squares above and on $B(n)^c$.
Finally, for each rook on $B(n)^c$, put a circle in those squares below and in
the column and on $B(n)^c$.
Then $\xi$ is
the number of squares with circles, minus the number of squares with both
circles and $x$'s, or $\# O - \# XO$ say.  Now directly from the board we have 
$$
\# O = \sum_{\sigma_i>i}n-\sigma_i+i
+\sum_{\sigma_i\le i} i-\sigma_i
$$
and
$$
\multline
\# XO = \# \{ 1\le i<j \le n; \, \sigma_i \le j < \sigma_j \} \\
+ \# \{ 1\le i<j \le n; \, \sigma_i > \sigma_j > j\} \\
+ \# \{ 1\le i<j \le n; \, j \ge \sigma_i > \sigma_j \}.
\endmultline
$$
The formula for $\# XO$ above equals $den(\sigma)$ [FoZe, p.33] and the 
formula for $\# O$ simplifies to $n\times \text{exc}(\sigma)$, so we get 
$$
\xi=n\times \text{exc}(\sigma) - \text{den}(\sigma)
$$ and the joint distribution of $(exc,den)$ with $(des,maj)$ follows
from (16) and the symmetry of the LHS of (13). 

%We can also perform an analysis identical to the one above except
%using $mat$ instead of $\xi$, yielding an explicit formula 
%which we list as a theorem.  It does not seem to be as easy to obtain
It is interesting to compare the statistic $den$ with the
following result, obtained by
performing the above analysis with $mat$ instead of $\xi$.  We list 
this in part $a$ of the theorem below.
 As far as the author is
aware, this is not equivalent to any known statistic.
\proclaim{Theorem 5}  

a)  For $\sigma_1\cdots \sigma_n \in S_n$, define
$$
\multline
\text{stat}(\sigma):=
\sum_{\sigma_i>i}\sigma_i-i
+\sum_{\sigma_i\le i} 1-\sigma_i  \\
+ \# \{ 1\le i<j \le n; \, \sigma_i > \sigma_j > j\}
+ \# \{ 1\le i<j \le n; \, \sigma_i \le j \text{ and }
\sigma_i < \sigma_j \}.
\endmultline
$$
Then $(exc,stat)$ is jointly distributed with $(des,maj)$.

b) For $\sigma \in M(\bold v)$, with $n=\sum_i v_i$, define
$$
\multline
\text{stat}x(\sigma):={n \choose 2} + \sum_{\sigma_i \le f_i(\bold v)} (
\sum_{i<j \atop \sigma_i > \sigma_j} 1 + \sum_{m<i \,  \,
\text{and} \, \, m\le v_1+\ldots 
+v_{\sigma_i-1} \atop
\sigma_m <\sigma_i}1) \\
+\sum_{\sigma_i > f_i(\bold v)}
\sum_{m<i \atop \sigma_m < \sigma_i}1 -
\sum_{\sigma_i > f_i(\bold v)}(i-1) -\sum_{\sigma_i \le f_i(\bold v)}
(n-i+\sum_{m<\sigma_i}v_m).
\endmultline
$$
Then $(exc,statx)$ is jointly distributed with $(des,maj)$.
\endproclaim
\noindent
$Proof:$  First we prove part $a$, where $\sigma \in S_n$.
For each rook in $C(\sigma)^T$, put an $x$ on the grid in all the
squares to the right and in the row.
For each rook off $B(n)^c$, put an $x$ in the squares below and in the
column, and also 
in the squares above and on $B(n)^c$.
Finally, for each rook on $B(n)^c$, put an $x$ in those squares above and in
the column.
Then $cross$ equals 
the number of rooks, plus the total number of $x$'s, minus the number of
squares with two $x$'s, or $\text{cross}=n+\#X-\#XX$ say.  Now
$$
n + \#X = 
{n+1 \choose 2} +  
\sum_{\sigma_i>i}n-\sigma_i+i
+\sum_{\sigma_i\le i} \sigma_i-1,
$$
and
$$
\multline
\# XX = \# \{ 1\le i<j \le n; \, \sigma_i > \sigma_j >j\} \\
+ \# \{ 1\le i<j \le n; \, \sigma_i \le j < \sigma_j \} \\
+ \# \{ 1\le i<j \le n; \,  \sigma_i < \sigma_j \le j \}.
\endmultline
$$
Since 
$$
\text{mat}(C(\sigma)^T,B(n)^c)=n\times \text{exc}(\sigma)+{n+1 \choose 2}
-\text{cross},
$$
the result follows.

%\input epsf
%\midinsert
%\def\epsfsize#1#2{.50#1}%make all figures 1/2 size
%$$\vbox{\epsffile{Figh15.eps}}$$
%\botcaption{Figure 15}  
%\endcaption
%\endinsert

It is desirable to have a multiset version of the explicit formula
from part $a$.  However, trying to mimic the above argument while using the
board $G_{\bold v}^c$ doesn't seem to lead to a nice formula.  Instead we
use $C(\sigma)$ and (19), and otherwise proceed exactly as in the
proof of part $a$.  The result is part $b$ above.
The details are left as an exercise to the interested
reader.
% Stat$x$ quickly reduces to $n^2-cross$.  
%Now place $x$'s in the squares of the grid exactly as in
%the proof of part a).  The 
%number of squares on or to the right of a rook is ${n+1 \choose 2}$.
%The number above a rook off $G_{\bold v}$ and on $G_{\bold v}$ is
%$\sum_{\sigma_i > f_i(\bold v)}(i-1)$.  The number below a rook off 
%$G_{\bold v}$ or above such a rook and on $G_{\bold v}$ is
%$\sum_{\sigma_i \le f_i(\bold v)} n-i+\sum_{m<\sigma_i}v_m$.  Hence
%$$
%\text{stat}x:={n \choose 2} 
%%+ \sum_{\sigma_i \le f_i(\bold v)} (
%%\sum_{i<j \atop \sigma_i < sigma_j} 1 + \sum_{i<j\le v_1+\ldots 
%v_{\sigma_i-1}}1) +\sum_{\sigma_i > f_i(\bold v)}
%\sum_{m<i \atop \sigma_m < \sigma_i}1 
%-
%\sum_{\sigma_i > f_i(\bold v)}(i-1) -\sum_{\sigma_i \le f_i(\bold v)}
%(n-i+\sum_{m<\sigma_i}v_m) + \#XX,
%$$
%where $\#XX$ is the number of squares with two $x$'s.  Now
%the number of squares with two $x$'s which are off $G_{\bold v}$ is
%$$
%\sum_{\sigma_i \le f_i(\bold v)}\sum_{i<j \atop \sigma_i>\sigma_j}1,
%$$
%and the number above a rook off $G_{\bold v}$ and on $G_{\bold v}$ is
%$$
%\sum_{\sigma_i \le f_i(\bold v)}\sum_{m<i \text{ and } \sigma_m<\sigma_i
%\atop m\le v_1+v_2+\ldots v_{\sigma_i-1}}.
%$$
%Finally, the number
%above a rook on $G_{\bold v}$ and on $G_{\bold v}$ is
%$$
%\sum_{\sigma_i > f_i(\bold v)}\sum_{m<i \atop \sigma_m<\sigma_i}1,
%$$
%and part b) of the theorem follows. 
$\qquad \blacksquare$

\head 5. Unimodality\endhead
In this section we show that for any admissible Ferrers board
$B$, $T_k(B)$ is symmetric and unimodal.  (A different proof that
$T_k(B)$ is symmetric can be found in [Dwo,p.52]).  For certain boards
we prove a stronger result.
\proclaim{Definition 8}  Let $f(q):=\sum_{j=M}^N a_jq^j$ be a
polynomial in $q$, where $a_M\ne 0$ and $a_N\ne 0$.  We call $f$
$symmetric$ if $a_{M+k}=a_{N-k}$ for $0\le k \le N-M$, and 
$unimodal$ if there exists $p$ such that $M\le p \le N$ and
$a_m\le a_{M+1}\le \cdots \le a_p\ge a_{p+1} \ge \cdots \ge a_N$.
Let $darga(f):=M+N$.  We say $f$ is $zsu(d)$ if $f$ is either
\vskip .1in
a) identically zero
\vskip .1in
\noindent
or
\vskip .1in
b) is $\in \Bbb N[q]$, and is
symmetric and unimodal with $darga(f)=d$.
\vskip .1in
\noindent
Note that the polynomial $q^s$ is $zsu(2s)$.
\endproclaim
\proclaim{Claim 1} If $f$ and $g$ are polynomials which are
both $zsu(d)$, then so is $f+g$. \endproclaim
\noindent
$Proof:$ Trivial. $\blacksquare$
\proclaim{Claim 2} If $f$ is $zsu(d)$ and $g$ is $zsu(e)$, then
$fg$ is $zsu(d+e)$. \endproclaim
\noindent
$Proof:$  (This proof is taken from [Zei]).  If either $f$ or $g$ is zero,
then so is $fg$.  If not, then $f$ can be written as a sum of ``atoms"
(terms of the form $q^{d-i}+q^{d-i+1}+\ldots +q^i$ for some $d/2\le i\le d$),
and $g$ equals the sum of atoms of the form $q^{e-j}+\ldots +q^j$.  The
product of two of these atoms is of the form
$$
q^{e+d-i-j}+2q^{e+d-i-j+1}+3q^{e+d-i-j+2}+\ldots +
3q^{i+j-2}+2q^{i+j-1}+q^{i+j},
$$
which is $zsu(d+e)$.  Summing over all products of atoms from $f$ and $g$,
and applying Claim 1 repeatedly proves the claim. $\blacksquare$

\proclaim{Definition 9}  Let $\delta$ be the linear operator such that
$\delta x^k:=[k]x^{k-1}$ for $k\in \Bbb Z$ ; for any formal power series $F(x)$,
$$
\delta F(x) = \frac{F(xq)-F(x)}{xq-x}.
$$
\endproclaim
\proclaim{Lemma 3}(This appears in [GaRe]).  For $0\le k \le n$,
$$
\delta \frac{ x^k }{ (1-x)(1-xq)\cdots (1-xq^n) } =
\frac{ [k]x^{k-1}+[n-k+1]x^k q }{ (1-x)(1-xq)\cdots (1-xq^{n+1}) }.
$$
\endproclaim
\proclaim{Definition 10}  
$$
\Phi (x;c_1,c_2,\ldots ,c_n):=\frac{\sum_{k=0}^n x^k 
T_{n-k}(B(c_1,\cdots ,c_n))}
{(1-x)(1-xq)\cdots (1-xq^n)},
$$
where $B(c_1,c_2,\ldots ,c_n)$ is the Ferrers board whose $i$th
column has height $c_i$.
\endproclaim
$\Phi$ satisfies the following useful identity [GaRe,p.259]
$$
\Phi (x;c_1,c_2,\ldots ,c_n)=\sum_{k=0}^{\infty}x^k
\prod_{i=1}^n [k+c_i-i+1]. \eqno(21)
$$
\proclaim{Theorem 6}  Let $B:=B(c_1,\ldots, c_n)$ be an 
admissible Ferrers board.  Then for $0\le k\le n$,
$$
\text{$T_k(B)$ is $zsu(N_k(B))$}, 
$$
where 
$$
N_k(B):= \text{Area}(B)
+n(n-k)-{n+1 \choose 2}. \eqno(22)
$$
\endproclaim
\noindent
$Proof:$  Throughout the proof, $B$ denotes the board
$B(c_1,c_2,\ldots ,c_n)$.  Our proof is a straight-forward refinement of the
proof in [GaRe,pp.258-263] that $T_k(B) \in \Bbb N [q]$.  First
we show that performing the board transformations RAISE, FLIP, and
ADD described below preserve property (22).  
\vskip .1in
\noindent
RAISE (this operation increases the height of each column by one;
it assumes $c_n\le n-1$).
\vskip .1in
\noindent Eq. (21) implies 
$$
\Phi (x;c_1+1,c_2+1,\ldots ,c_n+1)=
\Phi (x;c_1,c_2,\ldots ,c_n)/x,
$$
thus $T_k(B(c_1+1,c_2+1,\ldots ,c_n+1))=
T_{k-1}(B)$.  Hence if 
$T_{k-1}(B)$ is 
$zsu(\text{Area}(B)+n(n-k+1)-{n+1 \choose 2})$, then 
we have
$T_k(B(c_1+1,c_2+1,\ldots ,c_n+1))$ is  
$$
zsu(\text{Area}(B)+n+n(n-k)-{n+1 \choose 2})
$$
which is
$$
zsu(N_{k}(B(c_1+1,\ldots ,c_n+1))). \qquad \blacksquare
$$
\vskip .1in
\noindent
FLIP (this operation replaces $B$ by $B^{*}$, where $B^{*}$ is $B$
reflected about the cross diagonal, the same reflection utilized in
section 4).
\vskip .1in
\noindent Since $B^{*}$ has the same rook numbers as $B$ and hence
the same $q$-rook numbers as $B$,
$\Phi (x;B^*)=\Phi (x;B)$, and $T_k(B^*)=T_k(B)$.
Clearly $\text{Area}(B^*)=\text{Area}(B)$, and so $N_k(B^*)=N_k(B)$.
$\qquad \blacksquare$
\vskip .1in
\noindent
ADD (this operation adds a column of height zero to $B$).
\vskip .1in
Since 
$\Phi (x;0,c_1,c_2,\ldots ,c_n)=
x\delta x\Phi (x;c_1,c_2,\ldots ,c_n)$ [GaRe,p.260],
using Lemma 3 we get
$$
\Phi (x;0,c_1,c_2,\ldots ,c_n)=x\delta \sum_{k=1}^{n+1}
\frac{x^k T_{n-k+1}(B)}
{(1-x)(1-xq)\cdots (1-xq^n)}
$$
$$
=\frac{x}
{(1-x)(1-xq)\cdots (1-xq^{n+1})}\sum_{k=1}^{n+1} x^{k-1}[k]T_{n-k+1}(B)
+x^k[n-k+1]T_{n-k+1}(B)q^k,
$$
or
$$
\multline
\sum_{k=0}^{n+1}
\frac{x^k T_{n+1-k}(B(0,c_1,\ldots ,c_n))}
{(1-x)(1-xq)\cdots (1-xq^{n+1})}\\
=\sum_{k=1}^{n+1}
\frac{x^k ([k]T_{n+1-k}(B)+[n-k+2]q^{k-1}T_{n-k+2}(B))}
{(1-x)(1-xq)\cdots (1-xq^{n+1})}. \endmultline
$$
Comparing numerators we get
$T_{n+1}(B(0,c_1,\ldots ,c_n))=0$, and also (after 
replacing $k$ by $n+1-k$),
$$
T_k(B(0,c_1,\ldots ,c_n))=
[n+1-k]T_{k}(B)+[k+1]q^{n-k}T_{k+1}(B), \qquad 1\le k \le n+1. \eqno(23)
$$
Assuming $T_k(B)$ is $zsu(N_k(B))$, and also that $T_{k+1}(B)$ is
$zsu(N_{k+1}(B))$, both terms on the RHS of (23) have 
darga $N_{k}(B(0,c_1,\ldots ,c_n))$, and by Claim 1, the
LHS of (23) does also. $\qquad \blacksquare$

We now proceed with the proof of Theorem 6 by induction on
$\text{Area}(B)$.  If $\text{Area}(B)=0$, then $B$ is the trivial
board of width $n$ (i.e. $c_i=0$ for $1\le i \le n$).  It follows
from (3) that 
$$
\Phi (x;B) = \sum _{k=0}^n R_{n-k}(B)[k]! 
\frac{x^k}{(1-x)\cdots (1-xq^k)}
$$
$$
= [n]! 
\frac{x^n}{(1-x)\cdots (1-xq^n)},
$$
which implies 
$$
T_k(B)=
\cases [n]! \qquad \text{ if $k=0$} \\
0  \qquad \quad \text{    if $k>0$}
\endcases.
$$
Since $N_0(B)=n^2-{n+1 \choose 2}={n \choose 2}=\text{darga}([n]!)$, this 
shows Theorem 6 is true if $B$ is trivial.

The rest of the proof is precisely as in [GaRe].
Assume $\text{Area}(B)>0$, and that Theorem 6 is true for all
boards of smaller Area than $B$.  
We now show that $B$ can be
obtained from a board of smaller Area by a sequence of RAISE, FLIP,
or ADD operations, and Theorem 6 follows by induction.

\vskip .1in
\noindent
Case 1)  $1\le c_1 \le c_2 \le \cdots \le c_n$.  Applying
RAISE to $B(c_1-1,c_2-1,\ldots,c_n-1)$ results in $B$, and
$\text{Area}(B(c_1-1,\ldots,c_n-1))<\text{Area}(B)$.
\vskip .1in
\noindent
Case 2)  $0\le c_1\le c_2\le \cdots \le c_n=n$.  After performing
FLIP, $B^*$ falls under Case 1.
\vskip .1in
\noindent
Case 3) $0\le c_1\le c_2 \le \cdots \le c_n \le n-1$.  Let
$s:=\min \{i:c_i>0\}$ (if $s$ doesn't exist, $B$ is trivial).
\vskip .1in
Subcase a)  $n-s+1\ge c_n$.  Let $H:=B(c_s-1,c_{s+1}-1,\ldots ,
c_n-1)$, and note that $\text{Area}(H)<\text{Area}(B)$.  Since
$n-s+1\ge c_n$, performing RAISE to $H$ results in an
admissible board, and following this by $s-1$ ADD operations,
we end up with $B$.
\vskip .1in
Subcase b)  $n-s+1<c_n$.  After performing FLIP, $B^*$ falls
under Subcase a.
\noindent
This completes the proof of Theorem 6. $\qquad \blacksquare$

Let $B(h_1,d_1;h_2,d_2;\ldots ;h_t,d_t)$ denote the Ferrers 
board of Fig. 11.  Lemma 1 implies that 
$$
\frac{T_k(B)}{\prod_{i=1}^t [d_i]!} \in \Bbb N[q].
$$
\noindent
In certain cases we can show that $T_k(B)/\prod_i [d_i]!$ is
symmetric and unimodal; this is a stronger condition then
$T_k(B)$ being symmetric and unimodal by Claim 2.
\proclaim{Definition 11}  If $m<0$, extend the definition of
the $q$-binomial coefficient in the standard way;
$$
\left [ \matrix m \\ k \endmatrix \right ]:=
\frac{(1-q^m)(1-q^{m-1})\cdots (1-q^{m-k+1})}
{(1-q)(1-q^2)\cdots (1-q^k)}.
$$
We call $m$ the $numerator$ of the $q$-binomial coefficient.
Also given numbers $d_i,e_i$ and $h_i$ for $1\le i \le t$,
let $D_i,E_i$, and $H_i$ be abbreviations for the partial sums
$d_1+d_2+\ldots +d_i$,$e_1+\ldots +e_i$, and $h_1+\ldots +h_i$,
respectively, with $1\le i\le t$ and $D_0=E_0=H_0=0$.
\endproclaim
\proclaim{Claim 3}  Given integers $d_i$, $e_i$, and $h_i$,
with $0\le e_i \le d_i$, $d_i \in \Bbb P$, $h_i \in \Bbb N$
for $1\le i \le t$, let
$$
P(\bold e):=\prod_{i=1}^t
\left [ \matrix H_i-D_{i-1}+E_{i-1} \\ d_i-e_i \endmatrix \right ] 
\left [ \matrix D_i+D_{i-1}-H_i-E_{i-1} \\ e_i \endmatrix \right ].
$$
Set $d_0=0$ and assume that either

1) $d_{i-1}+d_i\ge h_i$ for $1\le i \le t$, 

\vskip .05in
\noindent
or

\vskip .05in
2) $D_i\ge H_i$ for $1\le i\le t$.

\noindent  Then if any of the numerators of the $q$-binomial
coefficients in the definition of $P(\bold e)$ are negative,
$P(\bold e)=0$.
\endproclaim
\noindent
$Proof:$  If $H_k-D_{k-1}+E_{k-1}<0$ for some $k$ with
$1\le k \le t$, choose $j$ so that for $i<j$, $H_i-D_{i-1}+E_{i-1}\ge 0$
and $H_j-D_{j-1}+E_{j-1}<0$.  Note that $j\ge 2$.  Now $H_j-D_{j-1}+E_{j-1}
<0$ implies $H_j-D_{j-2}+E_{j-2}<d_{j-1}-e_{j-1}$ which implies
$H_{j-1}-D_{j-2}+E_{j-2}<d_{j-1}-e_{j-1}$ which implies
$$
\left [ \matrix H_{j-1}-D_{j-2}+E_{j-2} \\ d_{j-1}-e_{j-1} \endmatrix \right ]
=0
$$
(since the numerator of this $q$-binomial coefficient is nonnegative by
definition of $j$).  This implies $P(\bold e)=0$.  Next assume we have a
$j$ for which $D_j+D_{j-1}-H_j-E_{j-1}<0$, but
$D_i+D_{i-1}-H_i-E_{i-1}\ge 0$ for $1\le i <j$.  If condition $2)$ is true,
this is impossible, for
$D_j+D_{j-1}-H_j-E_{j-1}\ge
D_j+D_{j-1}-H_j-D_{j-1}=D_j-H_j$.  So assume $1)$ holds.  Then
$D_j+D_{j-1}-H_j-E_{j-1}<0$ implies
$D_j+D_{j-1}-H_j-E_{j-2}<e_{j-1}$ which implies $d_j+d_{j-1}-h_j+
D_{j-1}+D_{j-2}-H_{j-1}-E_{j-2}<e_{j-1}$ which implies
$D_{j-1}+D_{j-2}-H_{j-1}-E_{j-2}<e_{j-1}$ which implies
$$
\left [ \matrix D_{j-1}+D_{j-2}-H_{j-1}-E_{j-2} \\ e_{j-1} 
\endmatrix \right ]=0,
$$
since the numerator is nonnegative by definition of $j$ but less than
the denominator.  $\qquad \blacksquare$
\proclaim{Claim 4}  Given integers $d_i$, $e_i$, and $h_i$ as in
Claim 3, let
$$
Q(s):=\prod_{i=1}^t 
\left [ \matrix s+H_i-D_{i-1} \\ d_{i} \endmatrix \right ].
$$
Then if $s+H_i-D_{i-1}<0$ for any $i$ satisfying $1\le i \le t$, then
$Q(s)=0$.
\endproclaim
\noindent
$Proof:$  Let $j$ be such that $s+H_j-D_{j-1}<0$, but $s+H_i-D_{i-1}
\ge 0$ for $1\le i <j$.  Note that $j \ge 2$ since $s+H_1-D_0=
s+h_1\ge 0$.  Now $s+H_j-D_{j-1}<0$ implies $s+H_j-D_{j-2}<d_{j-1}$
which implies $s+H_{j-1}-D_{j-2}<d_{j-1}$ which implies
$\left [ \matrix s+H_{j-1}-D_{j-2} \\ d_{j-1} \endmatrix \right ]=0$
since the numerator of this $q$-binomial coefficient is nonnegative.
$\qquad \blacksquare$

We have previously assumed that $B$ is an admissible board
($c_n \le n$) but in the next theorem we remove that restriction.  Note
that the definition of $R_k$ makes sense if $c_n > n$ as well, and for
such inadmissible boards we define $T_k$ via (3) (in general
these $T_k \notin \Bbb N[q]$).

\proclaim{Theorem 7}  Let $B=B(h_1,d_1;\ldots ;h_t,d_t)$ be the
Ferrers board of Fig. 11, where $H_t$ may be greater than $D_t$
($B$ inadmissible).  Set $L_k(B)=\text{Area}(B)+n(n-k)-
\sum_{i=1}^t D_id_i$.  Then $T_k(B)$ is either zero or
symmetric with $darga$ $L_k(B)$.  In addition, if either

1) $d_{i-1}+d_i\ge h_i$ for $1\le i \le t$, 

\vskip .05in
\noindent
or
\vskip .05in

2) $D_i\ge H_i$ for $1\le i\le t$,

\noindent then
$$
\frac{T_k(B)}{\prod_{i=1}^t [d_i]!} \text{ is }zsu(L_k(B)).
$$
\endproclaim
\noindent
$Proof:$  We require the following formulas
$$
\frac{T_{n-k}(B)}{\prod_{i=1}^t [d_i]!}=
\sum_{s=0}^{k} 
\left [ \matrix n+1 \\ k-s 
\endmatrix \right ] (-1)^{k-s} q^{k-s \choose 2}
\prod_{i=1}^t 
\left [ \matrix s+H_i-D_{i-1} \\ d_i 
\endmatrix \right ], \eqno(24)
$$
and
$$
\medmuskip0mu
T_{n-k}(B)=[d_t]! 
\sum_{k-d_t \le s \le k}T_{n-d_t-s}(B^{\prime})
\left [ \matrix H_{t}-n+d_t+s \\ d_{t}-k+s \endmatrix \right ]
\left [ \matrix 2n-d_t-H_t-s \\ k-s \endmatrix \right ]
q^{(k-s)(H_t+k-n)}, \eqno (25) 
$$
where $B^{\prime}=B(h_1,d_1;\ldots ;h_{t-1},d_{t-1})$ is obtained
by truncating the last $d_t$ columns of $B$.  The initial conditions
are given by $T_s(\emptyset)$ equals
$1$ if $s=n$ and zero otherwise, where $\emptyset$ denotes the
empty board with zero columns.

Eq. (24) is easily derived from (2), (3), and the
$q$-Vandermonde convolution [Hag1,p.98], [Dwo,p.39].
  Eq. (25) is Theorem
4.3.13 of [Hag1], and can also be obtained by setting $p=t,x=y=1$ in
Corollary 5.10 of [Hag2], where an inductive proof of the result
is given.
\proclaim{Lemma 4} 
With $B$ as above and $P(\bold e)$ as in
Claim 3,
$$
T_{n-k}(B)=\prod_{i=1}^t [d_i]! 
\sum_{e_1+e_2+\ldots +e_t=k \atop 0\le e_i \le d_i}
\prod_{i=1}^t P(\bold e)q^{e_i(H_i-D_i+E_i)}. \eqno (26)
$$
\endproclaim
\noindent
$Proof:$  By induction, the case $t=1$ following from (25).  
For $t>1$,
using (25) and
the inductive hypothesis we get
$$
\medmuskip0mu
T_{n-k}(B)=\prod_{i=1}^t [d_i]! \sum_{0 \le e_t \le k \atop
0\le e_t \le d_t} 
\left [ \matrix H_{t}-D_{t-1}+k-e_t \\ d_{t}-e_t \endmatrix \right ]
\left [ \matrix D_{t}+D_{t-1}-H_t-k+e_t \\ e_t \endmatrix \right ]
q^{e_t(H_t+k-D_t)} 
$$
$$
\times
\sum_{e_1+\ldots +e_{t-1}=k-e_t \atop
0\le e_i \le d_i}  \prod_{i=1}^{t-1}
\left [ \matrix H_{i}-D_{i-1}+E_{i-1} \\ d_{i}-e_i \endmatrix \right ]
\left [ \matrix D_{i}+D_{i-1}-H_i-E_{i-1} \\ e_i \endmatrix \right ]
q^{e_i(H_i-D_i+E_i)} 
$$
which equals the RHS of (26) since $k-e_t=E_{t-1}$.

We now proceed with the proof of Theorem 7.  We use the well-known
fact (see [GoOH],[Zei] for an amazing constructive proof) that 
for all $m\in \Bbb N$ and $k \in \Bbb N$, 
$\left [ \matrix m \\ k \endmatrix \right ]$ is $zsu(k(m-k))$.
Claim 4 implies that all the terms on the RHS of (24) are
polynomials.
  After a short calculation, we see that they are all
symmetric with $darga$ $L_k(B)$.
This proves the first part of the theorem, but unfortunately
the terms on the RHS of (24) alternate in sign and so we
cannot conclude that the LHS is unimodal.  However, if
condition 1) or 2) of Theorem 7 are satisfied, we can
apply Claim 3 and conclude all the terms
on the RHS of (26) are $\in \Bbb N[q]$. They are also all of $darga$ $L_k(B)$,
and so the second part of Theorem 7 
follows by Claim 1. $\blacksquare$
\proclaim{Corollary 3} For any vector $\bold v$ of 
nonnegative integers,
$$
\sum_{\pi \in M(\bold v) \atop \text{des}(\pi)=k} q^{\text{maj}(\pi)}
\text{ is }zsu(nk). 
$$
\endproclaim
\noindent 
$Proof:$ The board $G_{\bold v}$ satisfies condition 2) of Theorem 7,
and combining this with (18) we have
$$
\sum_{\pi \in M(\bold v) \atop \text{des}(\pi)=k} q^{\text{maj}(\pi)}
$$
is 
$$
zsu(\text{Area}(G_{\bold v}) +n(n-k)-\sum_{i=1}^tv_i(v_1+v_2+\ldots 
+v_i) +2nk -2\text{Area}(G_{\bold v})). \eqno(27)
$$  
Now 
$$
\text{Area}(G_{\bold v})=\sum_{i=1}^tv_i(v_1+v_2+\ldots +v_{i-1}),
$$
and since $(v_1+v_2+\ldots +v_t)^2=n^2$, (27) reduces to
$zsu(nk)$. $\qquad \blacksquare$.

\head 6. Final Comments\endhead

In [Hag1,p.130], the following more general form of Corollary 3 is derived 
$$
\sum_{\pi \in M(\bold v) \atop \text{$k$ $r$-descents}} q^{\text{$r$maj}(\pi)}
$$
is $zsu(nk+\sum_{i=1}^t v_i(v_{i-r+1}+\ldots +v_{i-1}))$,   
which involves the $(q-r)$ Simon Newcomb numbers introduced by
Rawlings [Raw].  The author hopes to describe connections
between these numbers and $q$-rook polynomials more fully elsewhere [Hag3].

Galovich and White have introduced a very general method of generating Mahonian
statistics, statistics they call ``splittable" [GaWh].  The author would like
to thank them for consultations regarding the statistic $mat(B)$, which together with
simple examples have led to the conclusion that $mat$ is not splittable,
at least not for all boards $B$.

For some time researchers have sought a $q$-analog of the theory of
permutations with restricted position.
 No positive
answer to this question has ever been found.  
Joni and Rota [JoRo] showed how the study of vector spaces over finite
fields with restricted bases is relevant to this problem. Later
Chen and Rota [ChRo] proved that if you require a $q$-analog 
 to have 
a certain interpretation in terms of automorphisms with
prescribed behavior, then a solution is possible only for a few types of 
boards.  There are interesting similarities between,
but no obvious overlap with, some of their
results and ours.

One can also try and develop a $q$-analog by finding a way of defining
$R_k$ for arbitrary boards (not just Ferrers boards) such that a
$q$-analog of (3) holds.  Perhaps the connection between
matrices over $\Bbb F_q$ of fixed rank and rook placements will shed
some light on this question.

\enddocument